\documentclass[preprint,12pt]{elsarticle}
\journal{Comput. Methods Appl. Mech. Eng.}
\usepackage{graphics,epsfig}
\usepackage{graphicx}
\usepackage{float}
\usepackage{amssymb,amstext,amsmath}
\usepackage{mathtools,physics,verbatim}
\usepackage{booktabs}
\usepackage{xspace}
\usepackage{lscape}
\usepackage{multirow}
\usepackage{lineno}
\usepackage[htt]{hyphenat}
\usepackage{xcolor}
\usepackage[toc,page]{appendix} 
\usepackage{hyperref}
\usepackage{cleveref}

\crefname{appsec}{Appendix}{Appendices}

\allowdisplaybreaks

\begin{document}

\begin{frontmatter}
\title{Asymptotically accurate and geometric locking-free finite element implementation of a refined shell theory}
\author{Khanh Chau Le$^{a,b}$\footnote{Corresponding author. Phone: +84 93 4152458, email: lekhanhchau@tdtu.edu.vn} \textnormal{and} Hoang-Giang Bui$^c$}
\address{$^a$Mechanics of Advanced Materials and Structures, Institute for Advanced Study in Technology, Ton Duc Thang University, Ho Chi Minh City, Vietnam
\\
$^b$Faculty of Civil Engineering, Ton Duc Thang University, Ho Chi Minh City, Vietnam
\\
$^c$Institute of Material Systems Modeling, Helmholtz-Zentrum Hereon, Geesthacht, Germany}
\begin{abstract} 
Accurate finite element analysis of refined shell theories is crucial but often hindered by membrane and shear locking effects. While various element-based locking-free techniques exist, this work addresses the problem at the theoretical level by utilizing results from asymptotic analysis. A formulation of a 2D refined shell theory incorporating transverse shear is developed using rescaled coordinates and angles of rotation, ensuring equal asymptotic orders of magnitude for extension, bending, and rotation measures and their respective stiffnesses. This novel approach, implemented via isogeometric analysis, is shown to be both asymptotically accurate relative to the underlying refined shell theory and inherently free from membrane and shear locking. Numerical simulations of semi-cylindrical shells show excellent agreement between the analytical solutions, 2D refined shell theory predictions, and 3D elasticity theory, validating the effectiveness and accuracy of the proposed formulation.
\end{abstract}


\end{frontmatter}


\section{Introduction}
Shells, thin-walled structures with curved surfaces, are ubiquitous in engineering, from civil and environmental to mechanical and aerospace applications. Understanding their structural behavior is crucial in designing safe, efficient and resilient structures. When a shell's thickness is significantly smaller than its characteristic radius of curvature and longitudinal size, its deformation can be effectively approximated using functions defined on 2D surface coordinates. Asymptotic analysis, using shell thickness as a model parameter, provides a hypothesis-free and systematic way to derive such 2D shell theories from 3D elasticity theory. This rigorous approach, based on the asymptotic analysis of 3D elasticity's strong \cite{gol2014theory} and weak formulation \cite{berdichevsky2009variational}, yields not only the classical shell theory \cite{sanders1959improved,koiter1960consistent} but also various refined shell theories (see \cite{gol2014theory,berdichevsky2009variational,berdichevsky1992effect}).

Accurate modeling of shells, especially concerning transverse shear effects, has driven the development of refined theories. Berdichevsky's pioneering work constructed an asymptotically exact refined shell theory using the variational-asymptotic method (VAM) \cite{berdichevsky1979variational,berdichevskii1980variational,berdichevsky2009variational}. This theory maintains the accuracy up to $O(h/R)$, $O(h/l)$, and $O(h^2/l^2)$, where $h$, $R$, and $l$ represent the shell thickness, curvature radius, and longitudinal deformation scale, respectively. While sharing asymptotic equivalence with Reissner's first-order shear deformation theory (FSDT) for plates \cite{reissner1945effect}, Berdichevsky’s formulation uniquely provides pointwise accuracy for both displacements and stresses, contrasting with Reissner’s integral-based accuracy. The necessity of refined theories is evident from the potential inaccuracies in displacement predictions by classical shell theory \cite{berdichevsky1992effect}. Subsequent efforts have expanded VAM applications to laminated structures \cite{sutyrin1997derivation,yu2005mathematical,yifeng2012avariational,yifeng2012variational}, often optimizing FSDT parameters for near-asymptotic correctness. Beyond linear elasticity, VAM has been applied to various nonlinear materials \cite{burela2012vam,burela2019asymptotically,bhadoria2024analytical,bhadoria2024asymptotically}. Recent investigations have explored asymptotically exact dimension reduction and error estimation for functionally graded plates \cite{le2023asymptotically} and beams \cite{le2025asymptotically}.
 
The mathematical complexity of 2D refined shell theories often precludes analytical solutions, making finite element implementation crucial for practical problems. However, membrane and shear locking can significantly reduce the reliability \cite{belytschko1985stress,bucalem1997finite,bletzinger2000unified}. Shear locking arises from the disparity between shear and bending stiffnesses, coupled with much smaller rotations from pure shear versus bending-induced curvature changes. This discrepancy, evident in the variational-asymptotic analysis of the energy functional \cite{berdichevsky2009variational,le1999vibrations}, leads to locking effect as shell thickness vanishes ($h\to 0$) with standard low-order elements. Membrane locking stems from the contrast between extension and bending stiffnesses and measures, causing numerical artifact due to multiplication of small and large quantities as $h\to 0$. As both membrane and shear locking occur in this limit, we commonly call them geometric locking.

While sophisticated methods exist to alleviate locking, they often compromise computational efficiency. Reduced and selective integration \cite{zienkiewicz1971reduced,onate1978techniques,hughes1978reduced}, using different rules for different energy components, is a popular technique but can introduce rank deficiency and spurious energy modes \cite{hughes2012finite}. Thus, alternative element-based techniques have been developed for improved accuracy and robustness. These include, among others, the modified shear strain method \cite{hughes1981finite}, hybrid/mixed formulations interpolating multiple fields \cite{saleeb1990hybrid,lavrenvcivc2021hybrid}, the assumed natural strains method \cite{simo1990class,caseiro2014assumed}, enhanced assumed strains methods \cite{buchter1994three}, and discrete strain gap approaches \cite{bletzinger2000unified,koschnick2005discrete,nguyen2013cell}, with some methods specifically adapted for isogeometric analysis \cite{caseiro2014assumed,oesterle2016shear}.

Building upon a recently developed, inherently locking-free and asymptotically accurate rescaled formulation for Berdichevsky's plate theory using $C^1$-isogeometric analysis \cite{le2024asymptotically}, this work extends the approach to refined shell theories. As established in the preceding discussion, standard finite element formulations of such theories suffer from geometric locking due to the inherent disparity in asymptotic scaling between bending and extension/shear contributions. Recognizing that this ill-conditioning is the fundamental barrier to achieving asymptotic accuracy numerically, our deliberation focused on addressing this issue directly at the variational level, prior to discretization. Inspired by the variational-asymptotic method's emphasis on correct scaling \cite{berdichevsky1979variational,le1999vibrations}, we introduce a novel rescaled formulation of Berdichevsky's refined shell theory \cite{berdichevsky1979variational}. This rescaling systematically balances the asymptotic orders of all kinematic measures and effective stiffnesses within the energy functional, yielding an intrinsically locking-free and asymptotically well-conditioned formulation. Our primary goals are: (i) To present this robust, rescaled shell formulation. (ii) To perform the first finite element verification of the asymptotic accuracy of Berdichevsky's refined shell theory, leveraging Isogeometric Analysis (IGA) for the requisite $C^1$-continuity of displacements and rotations \cite{berdichevsky1979variational,le2023asymptotically}. Numerical simulations of semi-cylindrical shells, compared against analytical and 3D elasticity solutions, validate the achievement of asymptotic accuracy with the proposed approach.
    
The paper is structured as follows: Section 2 outlines the 2D refined shell theory and its variational principles. Section 3 details the rescaled, locking-free variational formulation. Section 4 describes the finite element implementation. Section 5 presents numerical examples, including semi-cylindrical shells under internal pressure and a multi-component structure, and Section 6 concludes the paper. The Appendices provide the complete theoretical foundation for our method. The first Appendix details the derivation of the 2D rescaled variational formulation via an asymptotically consistent implicit B-matrix. The second Appendix then presents the weak formulation in matrix-vector form, yielding the explicit B-matrix.

\section{2D refined shell theory: Kinematics and variational principle}

\begin{figure}[htb]
    \centering
    \includegraphics[width=7cm]{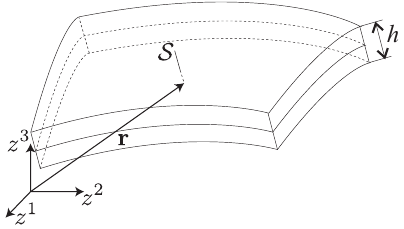}
    \caption{Schematic diagram of a shell segment}
    \label{fig:1}
\end{figure}

Consider a smooth, two-dimensional surface $\mathcal{S}$ within three-dimensional Euclidean space, bounded by a continuous closed curve $\partial \mathcal{S}$. We define a shell volume $\mathcal{V}$ by constructing line segments of length $h$, orthogonal to $\mathcal{S}$ at each point, with their midpoints on the surface. This definition is valid for sufficiently small shell thickness $h$, ensuring no segment intersections. Here, $\mathcal{S}$ represents the shell's middle surface, and $h$ its thickness (Fig.~\ref{fig:1}).

Mathematically, $\mathcal{S}$ is represented by the vector equation:
\begin{equation}
\vb{z}=\vb{r}(x^1,x^2),
\end{equation}
or, in component form:
\begin{equation}\label{surface}
z^i=r^i(x^\alpha ),\quad i=1,2,3;\quad \alpha=1,2,
\end{equation}
where $\vb{z}$ denotes the position vector in a 3D Cartesian system, and ${\bf r}(x^1,x^2)$ is a smooth vector function. Curvilinear coordinates $x^1$ and $x^2$ are chosen on $\mathcal{S}$, with units of length. We utilize Latin indices (1 to 3) for Cartesian coordinates and Greek indices (1 to 2) for surface coordinates.

The shell's geometric characteristics are described by the first and second fundamental forms:
\begin{align}
&a_{\alpha \beta }=\vb{t}_{\alpha }\vdot \vb{t}_{\beta }, \label{first}
\\
&b_{\alpha \beta }=-\vb{t}_{\alpha }\vdot \vb{n}_{,\beta 
}=\vb{t}_{\alpha ,\beta }\vdot \vb{n}, \label{second}
\end{align}
where $\vb{t}_{\alpha }=\vb{r}_{,\alpha }$ are tangent vectors, and $\vb{n}$ is the unit normal vector to $\mathcal{S}$, given by:
\begin{equation}
\vb{n}=\frac{\vb{t}_{1}\cross \vb{t}_{2}}{|\vb{t}_{1}| |\vb{t}_{2}|}.
\end{equation}
A comma preceding an index signifies partial differentiation.

Within 2D linear kinematics, shell deformation is defined by the displacement vector $\vb{u}(x^1,x^2)=u^i\vb{e}_i$, with $\vb{e}_i$ as Cartesian basis vectors. In the basis $\{ \vb{t}_{\alpha },\vb{n}\}$, displacement components are:
\begin{equation}
u_\alpha=\vb{t}_{\alpha }\vdot \vb{u},\quad u=\vb{n}\vdot \vb{u}.
\end{equation}
The 2D refined shell theory (2D-RST), considering transverse shear, introduces two additional degrees of freedom, $\varphi_\alpha (x^1,x^2)$, representing shear-induced rotation angles. The shear strain is a function of $\varphi_\alpha (x^1,x^2)$, with the explicit relationship given in Appendix A (Eq.~\eqref{strain2rs}). The shell's deformation is characterized by: (i) extension measures:
\begin{equation}
\label{extension}
\gamma_{\alpha \beta}=\vb{t}_{(\alpha }\vdot \vb{u}_{,\beta )}=u_{(\alpha ;\beta )}-b_{\alpha \beta }u,
\end{equation}
(ii) bending measures:
\begin{equation}
\label{bending}
\rho_{\alpha \beta}=(\vb{n}\vdot \vb{u}_{,(\alpha })_{;\beta)}+b^\lambda _{(\alpha }\varpi _{\beta )\lambda }-\varphi_{(\alpha ;\beta)}=u_{;\alpha \beta }+(u_\lambda b^\lambda _{(\alpha })_{;\beta )}+b^\lambda _{(\alpha }\varpi _{\beta )\lambda }-\varphi_{(\alpha ;\beta)},
\end{equation}
where
\begin{equation}
\varpi_{\alpha \beta }=\frac{1}{2}(\vb{t}_\beta \vdot \vb{u}_{,\alpha }-
\vb{t}_\alpha \vdot \vb{u}_{,\beta })=\frac{1}{2}(u_{\beta ,\alpha }-
u_{\alpha ,\beta }),
\end{equation}
and (iii) rotation angles $\varphi_\alpha$. When $\varphi_\alpha=0$, $\rho_{\alpha \beta}$ reduces to classical shell theory bending measures \cite{sanders1959improved,koiter1960consistent}. We employ the Einstein summation convention for repeated indices and use parentheses surrounding a pair of indices to denote symmetrization, as follows:
\begin{equation}
\label{parentheses}
u_{(\alpha ;\beta )}\equiv \frac{1}{2}(u_{\alpha ;\beta }+u_{\beta ;\alpha }), \quad b^\lambda _{(\alpha }\varpi _{\beta )\lambda }\equiv \frac{1}{2}(b^\lambda _{\alpha }\varpi _{\beta \lambda }+b^\lambda _{\beta }\varpi _{\alpha \lambda }).
\end{equation}
The metric tensor $a_{\alpha \beta }$ and its inverse $a^{\alpha \beta }$ are used to lower or raise  surface vector and tensor component indices. Covariant derivatives, indicated by a semicolon in indices, account for basis vector variations in surface tensors. For instance:
\begin{align}
&\nabla_\beta u_\alpha \equiv u_{\alpha ;\beta }=u_{\alpha ,\beta}-\Gamma^\lambda _{\alpha \beta }u_\lambda , \label{co-derivative}
\\
&\nabla_\beta u^\alpha \equiv u^\alpha _{;\beta }=u^\alpha _{,\beta}+\Gamma^\alpha _{\beta \lambda}u^\lambda , \label{contr-derivative}
\end{align}
where $\Gamma^\lambda _{\alpha \beta }$ are Christoffel's symbols:
\begin{equation}
\label{christoffel}
\Gamma ^\lambda _{\alpha \beta }=\frac{1}{2}a^{\lambda \mu }
(a_{\mu \alpha ,\beta }+a_{\mu \beta ,\alpha }
-a_{\alpha \beta ,\mu }).
\end{equation}
Note the difference in coordinate expressions for covariant derivatives of co- and contravariant vector components in \eqref{co-derivative}--\eqref{contr-derivative}. These definitions extend to tensor components \cite{le1999vibrations,duong2017a}.

The 2D-RST for linearly elastic, homogeneous shells is based on a variational principle \cite{berdichevsky2009variational,berdichevsky1979variational}. It states that the actual displacement field and rotation angles minimize the two-dimensional average energy functional
\begin{equation}\label{energy}
J[\vb{u},\varphi_\alpha ]=\int_{\mathcal{S} } \Phi(\gamma_{\alpha \beta},\rho_{\alpha \beta}, \varphi_\alpha ) \dd{\omega} - \mathcal{A}_{\rm{ext}}
\end{equation} 
among all admissible displacement fields and rotation angles that satisfy the kinematic boundary conditions. Here, $\dd \omega= \sqrt{a} \dd x^1 \dd x^2$ is the area element, and $a=\det a_{\alpha \beta }$. The two-dimensional energy density $\Phi(\gamma_{\alpha \beta},\rho_{\alpha \beta}, \varphi_\alpha )$ comprises three contributions, namely, (i) $\Phi_{\rm{cl}}(\gamma_{\alpha \beta},\rho_{\alpha \beta})$: energy density of classical shell theory, (ii) $\Phi_{\rm{gc}}(\gamma_{\alpha \beta},\rho_{\alpha \beta})$: geometric correction energy, (iii) $\Phi_{\rm{sc}}(\varphi_\alpha )$: shear correction energy. Thus:
\begin{equation}
\label{energydensity}
\Phi(\gamma_{\alpha \beta},\rho_{\alpha \beta},\varphi_\alpha ) =\Phi_{\rm{cl}}(\gamma_{\alpha \beta},\rho_{\alpha \beta})+\Phi_{\rm{gc}}(\gamma_{\alpha \beta},\rho_{\alpha \beta})+\Phi_{\rm{sc}}(\varphi_\alpha ) .
\end{equation}
These are explicitly given by:
\begin{align}
\Phi_{\rm{cl}}(\gamma_{\alpha \beta},\rho_{\alpha \beta})&=\mu h \Bigl[ \sigma (\gamma ^\alpha _\alpha)^2+\gamma_{\alpha \beta}\gamma^{\alpha \beta}\Bigr] +\frac{\mu h^3}{12} \Bigl[ \sigma (\rho^\alpha _\alpha)^2+\rho_{\alpha \beta}\rho^{\alpha \beta}\Bigr] , \label{cl1}
\\
\Phi_{\rm{gc}}(\gamma_{\alpha \beta},\rho_{\alpha \beta})&=-\frac{\mu h^3}{3} \Bigl[ \rho^{\alpha \beta}b^{\prime \lambda}_\alpha \gamma_{\beta \lambda }+\sigma \rho^{\alpha \beta}b_{\alpha \beta}\gamma^\lambda_\lambda+\frac{3}{5}\sigma \rho^\lambda_\lambda b^{\alpha \beta}\gamma_{\alpha \beta} \notag
\\
&+\sigma \Bigl( \frac{6}{5}\sigma-1\Bigr) \rho^\lambda_\lambda H \gamma^\mu _\mu \Bigr] , \label{gc1}
\\
\Phi_{\rm{sc}}(\varphi_\alpha )&=\frac{5}{12}\mu h a^{\alpha \beta} \varphi_\alpha \varphi_\beta , \label{sc1}
\end{align}
where $\lambda$ and $\mu$ are Lamé's constants, $\sigma= \frac{\lambda}{\lambda +2\mu}=\frac{\nu }{1-\nu}$, and $b^{\prime \lambda}_\alpha =b^{\lambda}_\alpha -H\delta^\lambda_\alpha$ is the deviator of the second quadratic form, with $H=b^\alpha_\alpha/2$ being the mean curvature. The work done by external loads, $\mathcal{A}_{\rm{ext}}$, is expressed as
\begin{multline}\label{work}
\mathcal{A}_{\rm{ext}}=\int_{\mathcal{S} } \Bigl[ (f_i-hHg_i)u^i+\frac{h}{2}g^\alpha (u_{,\alpha}+b^\lambda_\alpha u_\lambda )-\frac{\sigma h}{2}(g+\frac{1}{6}hf^\alpha_{;\alpha})\gamma^\beta_\beta
\\
+\frac{1}{10}\sigma h^2(f+\frac{1}{12}hg^\alpha_{;\alpha})\rho^\beta_\beta+\frac{1}{12}g^\alpha h\varphi_\alpha \Bigr] \dd{\omega},
\end{multline}
where
\begin{align}
&f_i = \tau_i |_{x_3=h/2} +\tau_i |_{x_3=-h/2},\quad f_\alpha=f_ir^i_{,\alpha},\quad f=f_in^i, \label{extforces}\\
&g_i = \tau_i |_{x_3=h/2} -\tau_i |_{x_3=-h/2},\quad g_\alpha=g_ir^i_{,\alpha},\quad g=g_in^i, \label{extforcesg}
\end{align}
with $\tau_i$ representing the Cartesian components of the traction vector. 

This 2D refined shell theory, by maintaining asymptotic accuracy up to $O(h/R)$ and $O(h^2/l^2)$, extends reliable analysis to shells of moderate thickness, where $R$ and $l$ denote the characteristic radius of curvature and longitudinal deformation length-scale, respectively \cite{berdichevsky1979variational,berdichevsky2009variational}.

\section{Rescaled formulation}
The finite element (FE) implementation of the 2D variational problem (Eq.~\eqref{energy}) presents challenges due to the differing orders of magnitude between bending and shear stiffnesses. Shear stiffness, significantly exceeding bending stiffness by two orders relative to shell thickness $h$, can cause numerical imbalance (shear locking) when multiplied by small rotation angles as $h$ decreases. A similar issue, membrane locking, arises from the contrast between extension and bending stiffnesses, and their respective measures.

To address these challenges, we introduce new unknown functions,
\begin{equation}
\psi_\alpha =-\vb{n}\vdot \vb{u}_{,\alpha }+\varphi_\alpha=-u_{,\alpha }-u_\lambda b^\lambda_\alpha +\varphi_\alpha ,
\end{equation} 
representing the total rotation angles of transverse fibers resulting from both bending and shear, such that
\begin{align}
&\rho_{\alpha \beta} = -\psi_{(\alpha;\beta )}+b^\lambda _{(\alpha }\varpi _{\beta )\lambda }, \label{rho}\\
&\varphi _\alpha= u_{,\alpha }+u_\lambda b^\lambda_\alpha +\psi_\alpha . \label{varphi}
\end{align} 
Consequently, the functional in Eq.~\eqref{energy} is reformulated with these unknowns as
\begin{equation}\label{energy2}
J[\vb{u},\psi _\alpha ]=\int_{\mathcal{S} } \Phi(\gamma_{\alpha \beta},\rho_{\alpha \beta},u_{,\alpha }+u_\lambda b^\lambda_\alpha +\psi_\alpha ) \dd{\omega} 
- \mathcal{A}_{\rm{ext}},
\end{equation} 
where function $\Phi$ of three arguments remains the same as in \eqref{energydensity}, while $\rho_{\alpha \beta}$ should be taken from \eqref{rho}. To avoid the membrane- and shear-locking effects, we aim to make the problem asymptotically independent of $h$, ensuring that extension, bending, and shear stiffnesses have comparable orders of magnitude. Following the approach first proposed in \cite{le2024asymptotically}, we introduce the rescaled coordinates:
\begin{equation}\label{scalingx}
\bar{z}^i=\frac{z^i}{h},\quad \bar{x}^\alpha =\frac{x^\alpha}{h},
\end{equation}
and express Eq.~\eqref{surface} as
\begin{equation}
\bar{z}^i=\bar{r}^i(\bar{x}^\alpha ),
\end{equation}
where $\bar{r}^i(\bar{x}^\alpha )=r^i(x^\alpha)/h$. This rescaling renders the left-hand sides of \eqref{scalingx} dimensionless. Importantly, the basis vectors and unit normal vector are unaffected, preserving the metric tensor and its inverse. Specifically,
\begin{equation}
\bar{\vb{r}}_{,\bar{\alpha}}=\pdv{\bar{\vb{r}}}{\bar{x}^\alpha}=\pdv{(\vb{r}/h)}{(x^\alpha/h)}=\vb{r}_{,\alpha},
\end{equation}
thus $\bar{\vb{n}}=\vb{n}$. 

Displacements are unchanged, but to reflect their dependence on the new base $\{\bar{\vb{r}}_{,\bar{\alpha}},\bar{\vb{n}}\}$ and argument $\bar{x}^\alpha$, they are denoted as:
\begin{equation}
\label{scalingu}
\bar{u}_{\bar{\alpha}}(\bar{x}^\alpha)=u_\alpha (x^\alpha),\quad \bar{u}(\bar{x}^\alpha)=u(x^\alpha).
\end{equation} 
The second fundamental form of the surface transforms as
\begin{equation}\label{scalingb}
\bar{b}_{\bar{\alpha} \bar{\beta}}=\bar{\vb{r}}_{,\bar{\alpha} \bar{\beta}}\vdot \bar{\vb{n}}=h \vb{r}_{,\alpha \beta }\vdot \vb{n}=hb_{\alpha \beta},
\end{equation} 
resulting in scaled principal radii of curvature $\bar{R}_1=R_1/h$ and $\bar{R}_2=R_2/h$. Partial derivatives with respect to $\bar{x}^\alpha$ are related to those with respect to $x^\alpha$ by
\begin{equation}
\pdv{\bar{x}^\alpha }=h\pdv{x^\alpha},\quad (.)_{,\bar{\alpha }}=h(.)_{,\alpha }.
\end{equation}
From Eq.~\eqref{christoffel}, the scaled Christoffel symbols are related to the original ones by
\begin{equation}\label{rescaledgamma}
\bar{\Gamma }^{\bar{\lambda }}_{\bar{\alpha }\bar{\beta }}=\frac{1}{2}\bar{a}^{\bar{\lambda} \bar{\mu }}
(\bar{a}_{\bar{\mu }\bar{\alpha },\bar{\beta }}+\bar{a}_{\bar{\mu }\bar{\beta },\bar{\alpha }}-\bar{a}_{\bar{\alpha }\bar{\beta },\bar{\mu }})=\frac{1}{2}ha^{\lambda \mu }
(a_{\mu \alpha ,\beta }+a_{\mu \beta ,\alpha }
-a_{\alpha \beta ,\mu })=h\Gamma ^\lambda _{\alpha \beta }.
\end{equation}
Therefore, covariant derivatives with respect to $\bar{x}^\alpha$ are also related to those with respect to $x^\alpha$ by
\begin{equation}
\bar{\nabla}_{\bar{\alpha }}=h\nabla_{\alpha},\quad (.)_{;\bar{\alpha }}=h(.)_{;\alpha }.
\end{equation}
For rotation angles, we introduce
\begin{equation}
\label{scalingpsi}
\bar{\psi}_{\bar{\alpha}}=h\psi_{\alpha},
\end{equation}
with dimensions of length, representing longitudinal displacements of the shell's positive face \cite{berdichevsky1979variational}.

Applying the scaling rules \eqref{scalingx}, \eqref{scalingu}, and \eqref{scalingpsi}, we obtain the following relationships between original and rescaled quantities:
\begin{align}
&\gamma_{\alpha \beta }=\frac{1}{h}\bar{\gamma}_{\bar{\alpha }\bar{\beta }}, \quad 
\rho_{\alpha \beta}=\frac{1}{h^2}\bar{\rho}_{\bar{\alpha }\bar{\beta}}, \label{scaling}
\\
&u_{,\alpha }+u_\lambda b^\lambda_\alpha +\psi_\alpha =\frac{1}{h}(\bar{u}_{,\bar{\alpha} } +\bar{u}_{\bar{\lambda }}\bar{b}^{\bar{\lambda}}_{\bar{\alpha}} +\bar{\psi}_{\bar{\alpha}}), \label{scalingvarphi}
\\
&\dd{\omega}=\sqrt{a}\dd{x^1}\dd{x^2}=h^2\sqrt{\bar{a}}\dd{\bar{x}^1}\dd{\bar{x}^2} =h^2\dd{\bar{\omega}}, \label{scalingomega}
\end{align}
where
\begin{align}
&\bar{\gamma}_{\bar{\alpha }\bar{\beta }}=\bar{u}_{(\bar{\alpha };\bar{\beta })}-\bar{b}_{\bar{\alpha }\bar{\beta }}\bar{u}, \label{scalinggamma}
\\
&\bar{\rho}_{\bar{\alpha }\bar{\beta}}= -\bar{\psi}_{(\bar{\alpha};\bar{\beta })}+\bar{b}^{\bar{\lambda }}_{(\bar{\alpha }}\bar{\varpi }_{\bar{\beta })\bar{\lambda }}, \label{scalingrho}
\end{align}
and
\begin{equation}
\bar{\varpi }_{\bar{\alpha }\bar{\beta }}=\frac{1}{2}(\bar{u}_{\bar{\beta },\bar{\alpha }}-\bar{u}_{\bar{\alpha },\bar{\beta }}).
\end{equation}

Substituting Eqs.~\eqref{scaling}--\eqref{scalingomega} into the energy functional \eqref{energy2} gives a form with rescaled quantities:
\begin{equation}\label{rescaledfunct}
J[\bar{\vb{u}},\bar{\psi }_{\bar{\alpha}} ]=\mu h \int_{\bar{\mathcal{S} }} \bar{\Phi}(\bar{\gamma}_{\bar{\alpha }\bar{\beta }},\bar{\rho}_{\bar{\alpha }\bar{\beta}}, \bar{u}_{,\bar{\alpha} } +\bar{u}_{\bar{\lambda }}\bar{b}^{\bar{\lambda}}_{\bar{\alpha}} +\bar{\psi}_{\bar{\alpha}}) \dd{\bar{\omega}} - \mathcal{A}_{\rm{ext}} .
\end{equation} 
Here $\bar{\mathcal{S}}=\{ (\bar{x}_1,\bar{x}_2) \, | \, (x_1,x_2)\in \mathcal{S}\}$ denotes the rescaled 2D domain, while the rescaled energy density, $\bar{\Phi}(\bar{\gamma}_{\bar{\alpha }\bar{\beta }},\bar{\rho}_{\bar{\alpha }\bar{\beta}}, \bar{u}_{,\bar{\alpha} } +\bar{u}_{\bar{\lambda }}\bar{b}^{\bar{\lambda}}_{\bar{\alpha}} +\bar{\psi}_{\bar{\alpha}})$, is the sum of three contributions 
\begin{align}
\bar{\Phi}_{\rm{cl}}(\bar{\gamma}_{\bar{\alpha }\bar{\beta }},\bar{\rho}_{\bar{\alpha }\bar{\beta}})&= \sigma (\bar{\gamma }^{\bar{\alpha}}_ {\bar{\alpha}})^2+\bar{\gamma}_{\bar{\alpha }\bar{\beta }}\bar{\gamma}^{\bar{\alpha }\bar{\beta }} +\frac{1}{12} \Bigl[ \sigma (\bar{\rho }^{\bar{\alpha}}_{\bar{\alpha}})^2+\bar{\rho}_{\bar{\alpha }\bar{\beta }}\bar{\rho}^{\bar{\alpha }\bar{\beta }}\Bigr] ,\label{cl3}
\\
\bar{\Phi}_{\rm{gc}}(\bar{\gamma}_{\bar{\alpha }\bar{\beta }},\bar{\rho}_{\bar{\alpha }\bar{\beta}})&=-\frac{1}{3} \Bigl[ \bar{\rho}^{\bar{\alpha }\bar{\beta }}\bar{b}^{\prime \bar{\lambda}}_{\bar{\alpha}} \bar{\gamma}_{\bar{\beta }\bar{\lambda }}+\sigma \bar{\rho}^{\bar{\alpha }\bar{\beta }}\bar{b}_{\bar{\alpha }\bar{\beta }}\bar{\gamma}^{\bar{\lambda}}_{\bar{\lambda}}+\frac{3}{5}\sigma \bar{\rho}^{\bar{\lambda}}_{\bar{\lambda}} \bar{b}^{\bar{\alpha }\bar{\beta }}\bar{\gamma}_{\bar{\alpha }\bar{\beta }} \notag
\\
&+\sigma \Bigl( \frac{6}{5}\sigma-1\Bigr) \bar{\rho}^{\bar{\lambda}}_{\bar{\lambda}} \bar{H} \bar{\gamma}^{\bar{\mu}}_{\bar{\mu}} \Bigr], \label{gc3}
\\
\bar{\Phi}_{\rm{sc}}(\bar{u}_{,\bar{\alpha} } +\bar{u}_{\bar{\lambda }}\bar{b}^{\bar{\lambda}}_{\bar{\alpha}} +\bar{\psi}_{\bar{\alpha}})&= \frac{5}{12}\bar{a}^{\bar{\alpha }\bar{\beta }}(\bar{u}_{,\bar{\alpha} } +\bar{u}_{\bar{\lambda }}\bar{b}^{\bar{\lambda}}_{\bar{\alpha}} +\bar{\psi}_{\bar{\alpha}})(\bar{u}_{,\bar{\beta} } +\bar{u}_{\bar{\mu}}\bar{b}^{\bar{\mu}}_{\bar{\beta}} +\bar{\psi}_{\bar{\beta}}) \label{sc3}
\end{align}
similar to \eqref{energydensity}. The work of external forces becomes
\begin{multline}
\label{rescaledwork}
\mathcal{A}_{\rm{ext}}=h^2\int_{\mathcal{S} } \Bigl[ (f_i-\bar{H}g_i)\bar{u}^i+\frac{1}{2}g^{\bar{\alpha}} (\bar{u}_{,\bar{\alpha} } +\bar{u}_{\bar{\lambda }}\bar{b}^{\bar{\lambda}}_{\bar{\alpha}})-\frac{\sigma }{2}(g+\frac{1}{6}f^{\bar{\alpha}}_{;\bar{\alpha}})\bar{\gamma}^{\bar{\beta}}_{\bar{\beta}}
\\
+\frac{1}{10}\sigma (f+\frac{1}{12}g^{\bar{\alpha}}_{;\bar{\alpha}})\bar{\rho}^{\bar{\beta}}_{\bar{\beta}}+\frac{1}{12}g^{\bar{\alpha}} (\bar{u}_{,\bar{\alpha} } +\bar{u}_{\bar{\lambda }}\bar{b}^{\bar{\lambda}}_{\bar{\alpha}} +\bar{\psi}_{\bar{\alpha}}) \Bigr] \dd{\bar{\omega}} .
\end{multline}

Since scaling a functional by a constant does not alter its minimizer, we further simplify \eqref{rescaledfunct} by dividing by $\mu h$. The minimization problem becomes
\begin{equation}\label{eq:1}
\bar{J}[\bar{\vb{u}},\bar{\psi }_{\bar{\alpha }} ]=\int_{\bar{\mathcal{S} }} \bar{\Phi}(\bar{\gamma}_{\bar{\alpha }\bar{\beta }},\bar{\rho}_{\bar{\alpha }\bar{\beta}}, \bar{u}_{,\bar{\alpha} } +\bar{u}_{\bar{\lambda }}\bar{b}^{\bar{\lambda}}_{\bar{\alpha}} +\bar{\psi}_{\bar{\alpha}}) \dd{\bar{\omega}} 
- \bar{\mathcal{A}}_{\rm{ext}} \rightarrow \min_{\bar{u}_{\bar{\alpha}},\bar{u},\bar{\psi}_{\bar{\alpha }}},
\end{equation}
where
\begin{multline}
\bar{\mathcal{A}}_{\rm{ext}} =\int_{\mathcal{S} } \Bigl[ (\bar{f}_i-\bar{H}\bar{g}_i)\bar{u}^i+\frac{1}{2}\bar{g}^{\bar{\alpha}} (\bar{u}_{,\bar{\alpha} } +\bar{u}_{\bar{\lambda }}\bar{b}^{\bar{\lambda}}_{\bar{\alpha}})-\frac{\sigma }{2}(\bar{g}+\frac{1}{6}\bar{f}^{\bar{\alpha}}_{;\bar{\alpha}})\bar{\gamma}^{\bar{\beta}}_{\bar{\beta}}
\\
+\frac{1}{10}\sigma (\bar{f}+\frac{1}{12}\bar{g}^{\bar{\alpha}}_{;\bar{\alpha}})\bar{\rho}^{\bar{\beta}}_{\bar{\beta}}+\frac{1}{12}\bar{g}^{\bar{\alpha}} (\bar{u}_{,\bar{\alpha} } +\bar{u}_{\bar{\lambda }}\bar{b}^{\bar{\lambda}}_{\bar{\alpha}} +\bar{\psi}_{\bar{\alpha}}) \Bigr] \dd{\bar{\omega}},
\end{multline}
and
\begin{equation} 
\label{forcesc}
\bar{f}_i=\frac{hf_i}{\mu },\quad \bar{g}_i=\frac{hg_i}{\mu },\quad \bar{f}=\frac{hf}{\mu },\quad \bar{g}=\frac{hg}{\mu },\quad \bar{f}^{\bar{\alpha}}=\frac{hf^\alpha}{\mu },\quad \bar{g}^{\bar{\alpha}}=\frac{hg^\alpha}{\mu }.
\end{equation}
Consequently, both $\bar{f}_i$ and $\bar{g}_i$ scale with shell thickness $h$ and strain $\varepsilon=\max \{ \bar{f}_i/\mu, \bar{g}_i/\mu \}$. Since all three stiffnesses (extension, bending, and shear) in the rescaled functional \eqref{eq:1} are $O(1)$, the minimizer is $O(\bar{f}_i$,$\bar{g}_i)$ times a function depending on the characteristic size of $\bar{\mathcal{S}}$. Returning to original functions, extension measures $\gamma_{\alpha \beta}$ and rotation angles $\psi_\alpha$ are $O(\varepsilon)$ (with $\varphi_\alpha$ smaller \cite{le1999vibrations,berdichevsky2009variational}), while bending measures $\rho_{\alpha \beta}$ are $O(\varepsilon/h)$. Note that the rescaled energy contributions \eqref{cl3}--\eqref{sc3} still contain terms dependent on $h$ through the rescaled second quadratic form $\bar{b}_{\bar{\alpha }\bar{\beta }}$. However, $\bar{b}_{\bar{\alpha }\bar{\beta }}\to 0$ in the limit $h\to 0$, so the energy functional exhibits regular asymptotic behavior, ensuring that the associated variational problem defined in Eq.~\eqref{eq:1} becomes asymptotically well-conditioned. The rescaled problem \eqref{eq:1} avoids membrane and shear locking.

To solve problem \eqref{eq:1}, boundary conditions must be specified. If a portion of the shell's edge, $\bar{\partial}_k$, is clamped, then we require the admissible functions to satisfy the following kinematic conditions:
\begin{equation}\label{kinbc}
\bar{u}_{\bar{\alpha }} =0,\quad \bar{u}=0, \quad \bar{\psi }_{\bar{\alpha }} =0 \quad \text{at $\bar{\partial }_k$}.
\end{equation}
For a simply supported edge $\bar{\partial}_{ss}$, kinematical conditions $\bar{u}_{\bar{\alpha }} =0$ and $\bar{u}=0$ are enforced, while $\bar{\psi }_{\bar{\alpha }}$ are allowed to vary freely. On a free edge $\bar{\partial}_s$, no constraints are imposed on $\bar{u}_{\bar{\alpha }}$, $\bar{u}$ and $\bar{\psi }_{\bar{\alpha }}$. Beyond these typical cases, Section 5 explores additional boundary conditions to validate our FE implementation.

After obtaining the solution, the theory's asymptotic accuracy is assessed by comparing the true average displacement with 3D elasticity results. According to \cite{berdichevsky1979variational},  
$\bar{u}_{\bar{\alpha}}$ represent the true average tangential displacements. The true average normal displacement is obtained by correcting the normal displacement:
\begin{equation}\label{trueu}
\check{u}=u+\frac{h^2\sigma }{60}a^{\alpha \beta }\rho_{\alpha \beta}=\bar{u}+\frac{\sigma }{60}\bar{a}^{\bar{\alpha} \bar{\beta }}(-\bar{\psi}_{(\bar{\alpha };\bar{\beta})}+\bar{b}^{\bar{\lambda }}_{(\bar{\alpha }}\bar{\varpi }_{\bar{\beta })\bar{\lambda }}).
\end{equation}
Since $\varphi_\alpha$ also depends on first derivatives of $u$, we seek solutions where  
$\bar{u}_{\bar{\alpha}}$, $\bar{u}$, and $\bar{\psi}_{\bar{\alpha}}$ are $C^1$-functions for accurate computations and comparisons. Asymptotic accuracy is confirmed by agreement up to $O(h/R)$ and $O(h^2/l^2)$ between functions $\bar{u}_{\bar{\alpha}}$, $\check{u}$, and $\bar{\psi}_{\bar{\alpha}}$  found by the 2D refined shell theory with
\begin{equation}\label{eq:ex2_displacement_rotation}
r^i_{,\alpha }\langle w_i (x^\alpha ,x)\rangle ,\quad n^i \langle w_i (x^\alpha ,x)\rangle , \quad \text{and} \quad r^i_{,\alpha }\langle w_i  (x^\alpha ,x) x \rangle /(h^2/12) ,
\end{equation}
where $\langle .\rangle \equiv \frac{1}{h}\int_{-h/2}^{h/2} . \dd{x}$ denotes the averaging over the shell's thickness and $w_i(x^\alpha,x)$ are the displacements computed by the 3D exact theory of elasticity.  

\section{Finite element implementation}

\subsection{Weak and strong formulations} \label{sec:weak_strong}
In this Section and the theoretical portion of the subsequent Section, we will use rescaled coordinates and quantities exclusively. For brevity, overbars will be omitted.

Taking the first variation of the functional (Eq.~\eqref{eq:1}), a necessary condition for its minimizer is that the virtual work of internal forces equals the virtual work of external forces:
\begin{equation}
\label{var1}
\int_{\mathcal{S}} \Bigl[ n^{\alpha \beta } \var{\gamma }_{\alpha \beta }+m^{\alpha \beta } \var{\rho }_{\alpha \beta }+q^\alpha (\var{u}_{,\alpha }+b^\lambda _\alpha \var{u}_\lambda +\var{\psi }_\alpha )\Bigr] \dd{\omega}=\var{\mathcal{A}}_{\rm{ext}},
\end{equation}
where the membrane forces $n^{\alpha \beta }$, bending moments $m^{\alpha \beta}$, and shear forces $q^\alpha$ are dual to $\gamma_{\alpha \beta}$, $\rho_{\alpha \beta}$ and $\varphi_\alpha$, respectively. These are defined as:
\begin{align}
n^{\alpha \beta }=\pdv{\Phi}{\gamma_{\alpha \beta}}&=2(\sigma \gamma^\lambda _\lambda a^{\alpha \beta }+\gamma^{\alpha \beta})-\frac{1}{3}\Bigl[ \rho^{(\alpha \lambda }b^{\prime \beta )}_\lambda+\sigma a^{\alpha \beta }\Bigl( b_{\mu \nu}\rho^{\mu \nu} \notag
\\
&+\Bigl(\frac{6}{5}\sigma -1\Bigr)H\rho^\lambda_\lambda \Bigr)+\frac{3}{5}\sigma \rho^\lambda _\lambda b^{\alpha \beta }\Bigr], \label{mforce}
\\
m^{\alpha \beta }=\pdv{\Phi}{\rho_{\alpha \beta}}&=\frac{1}{6}(\sigma \rho^\lambda _\lambda a^{\alpha \beta }+\rho^{\alpha \beta})-\frac{1}{3}\Bigl[ \gamma^{(\alpha \lambda }b^{\prime \beta )}_\lambda+\sigma a^{\alpha \beta }\Bigl( \frac{3}{5}b_{\mu \nu}\gamma^{\mu \nu} \notag
\\
&+\Bigl(\frac{6}{5}\sigma -1\Bigr)H\gamma^\lambda_\lambda \Bigr)+\sigma \gamma^\lambda _\lambda b^{\alpha \beta }\Bigr], \label{moments}
\\
q^\alpha =\pdv{\Phi}{\varphi_{\alpha }}&=\frac{5}{6}a^{\alpha \beta }(u_{,\beta }+b^\lambda _\beta u_\lambda +\psi _\beta ). \label{qforce}
\end{align}
The virtual work of external forces is:
\begin{multline}\label{virtualwork}
\var{\mathcal{A}}_{\rm{ext}}=\int_{\mathcal{S}} \Bigl[ (f^\alpha -Hg^\alpha)\var{u}_\alpha +(f-Hg)\var{u}+\frac{1}{2}g^\alpha (\var{u}_{,\alpha }+b^\lambda _\alpha \var{u}_\lambda )
\\
-\frac{\sigma}{2}(g+\frac{1}{6}f^\lambda_{;\lambda })a^{\alpha \beta}(\var{u}_{(\alpha ;\beta )}-b_{\alpha \beta} \var{u}) 
-\frac{\sigma}{10}(f+\frac{1}{12}g^\lambda_{;\lambda })a^{\alpha \beta}\var{\psi}_{\alpha ;\beta }
\\
+\frac{1}{12}g^\alpha (\var{u}_{,\alpha } 
+b^\lambda _\alpha \var{u}_\lambda +\var{\psi }_\alpha )\Bigr] \dd{\omega}.
\end{multline}

Introducing the notations
\begin{multline}\label{virtualextension}
\var{W}^{m}=\int_{\mathcal{S}} n^{\alpha \beta } \var{\gamma }_{\alpha \beta }\dd{\omega}= \int_{\mathcal{S}} \dd{\omega}\Bigl\{ 2(\sigma \gamma^\lambda _\lambda a^{\alpha \beta }+\gamma^{\alpha \beta})-\frac{1}{3}\Bigl[ \rho^{(\alpha \lambda }b^{\prime \beta )}_\lambda 
\\
+\sigma a^{\alpha \beta }\Bigl( b_{\mu \nu}\rho^{\mu \nu}+\Bigl(\frac{6}{5}\sigma -1\Bigr)H\rho^\lambda_\lambda \Bigr)+\frac{3}{5}\sigma \rho^\lambda _\lambda b^{\alpha \beta }\Bigr] \Bigr\} (\var{u}_{(\alpha ;\beta )}-b_{\alpha \beta} \var{u}) \dd{\omega},
\end{multline}
\begin{multline}\label{virtualbending}
\var{W}^{b}=\int_{\mathcal{S}} m^{\alpha \beta } \var{\rho }_{\alpha \beta }\dd{\omega}= \int_{\mathcal{S}} \dd{\omega}\Bigl\{ \frac{1}{6}(\sigma \rho^\lambda _\lambda a^{\alpha \beta }+\rho^{\alpha \beta})-\frac{1}{3}\Bigl[ \gamma^{(\alpha \lambda }b^{\prime \beta )}_\lambda 
\\
+\sigma a^{\alpha \beta }\Bigl( \frac{3}{5}b_{\mu \nu}\gamma^{\mu \nu}+\Bigl(\frac{6}{5}\sigma -1\Bigr)H\gamma^\lambda_\lambda \Bigr)+\sigma \gamma^\lambda _\lambda b^{\alpha \beta }\Bigr] \Bigr\} (-\var{\psi }_{(\alpha ;\beta )}+b^\lambda _{(\alpha } \var{\varpi }_{\beta )\lambda }) \dd{\omega},
\end{multline}
and
\begin{equation}
\label{virtualshear}
\var{W}^{s}=\int_{\mathcal{S}} q^\alpha \var{\varphi}_\alpha \dd{\omega}=\int_{\mathcal{S}} \frac{5}{6}a^{\alpha \beta }(u_{,\beta}+b^\lambda _\beta u_\lambda +\psi _\beta )(\var{u}_{,\alpha }+b^\mu _\alpha \var{u}_\mu +\var{\psi }_\alpha ) \dd{\omega},
\end{equation}
for the virtual work of membrane forces, bending moments, and shear forces, respectively, Eq.~\eqref{var1} becomes:
\begin{equation}
\label{var2}
\var{W}^{m}+\var{W}^b+\var{W}^s=\var{\mathcal{A}}_{\rm{ext}}.
\end{equation}
Let $\mathcal{K}=\{ (v,v_\alpha, \chi_\alpha )\, | \, (v,v_\alpha,\chi_\alpha )|_{\partial_k}=0\}$ be the space of kinematically admissible functions. The weak formulation of the refined shell theory states that for given $f,f^\alpha,g,g^\alpha$ find $(u,u_\alpha,\psi_\alpha )\in \mathcal{K}$ such that Eq.~\eqref{var2} is satisfied for all $(\var{u},\var{u}_\alpha,\var{\psi}_\alpha )\in \mathcal{K}$. Looking at the integrals \eqref{virtualextension}--\eqref{virtualshear} we see that, for this weak formulation to make sense, the admissible functions should belong at least to the Sobolev's space of square integrable functions with the square integrable first derivatives, $H^1(\mathcal{S})$. However, since the desired asymptotic accuracy of FSDT may require higher smoothness of $u$, $u_\alpha$, and $\psi_\alpha$, the continuity assumption in $\mathcal{K}$ still remains unspecified.

For completeness, we also present the strong formulation of the problem. Assuming $u$, $u_\alpha$, and $\psi_\alpha$ are twice differentiable, we integrate Eq.~\eqref{var1} by parts, yielding:
\begin{multline}
\label{var3}
\int_{\mathcal{S}} \Bigl[ -n^{\alpha \beta }_{;\beta }\var{u}_\alpha -b_{\alpha \beta} n^{\alpha \beta }\var{u}+m^{\alpha \beta }_{;\beta } \var{\psi }_{\alpha }- (m^{\lambda [ \alpha}b^{\beta ] }_\lambda )_{;\beta }\var{u}_\alpha -q^\alpha _{;\alpha }\var{u}+b^\alpha _\beta q^\beta \var{u}_\alpha 
\\
+q^\alpha \var{\psi }_\alpha )\Bigr] \dd{\omega} +\int_{\partial_s} \Bigl[ n^{\alpha \beta }\nu_{\beta }\var{u}_\alpha -m^{\alpha \beta }\nu_{\beta } \var{\psi }_{\alpha } +m^{\lambda [ \alpha}b^{\beta ] }_\lambda \nu_{\beta }\var{u}_\alpha +q^\alpha \nu_{\alpha }\var{u} \Bigr] \dd{s}
\\
=\int_{\mathcal{S}} \Bigl[ (f^\alpha -Hg^\alpha )\var{u}_\alpha +(f-Hg)\var{u}-\frac{1}{2}g^\alpha _{;\alpha }\var{u}+\frac{1}{2}g^\beta b^\alpha _\beta \var{u}_\alpha  
\\
+\frac{\sigma}{2}(g+\frac{1}{6}f^\lambda _{;\lambda })_{;\beta }a^{\alpha \beta}\var{u}_\alpha+\sigma H(g+\frac{1}{6}f^\lambda _{;\lambda })\var{u}  +\frac{\sigma}{10}(f+\frac{1}{12}g^\lambda _{;\lambda })_{;\beta }a^{\alpha \beta}\var{\psi }_\alpha 
\\
-\frac{1}{12}g^\alpha_{;\alpha}\var{u}+\frac{1}{12}g^\beta b^\alpha_\beta \var{u}_\alpha +\frac{1}{12}g^\alpha \var{\psi}_\alpha \Bigr] \dd{\omega} +\int_{\partial_s} \Big[ \frac{1}{2}g^\alpha \nu _{\alpha }\var{u}
\\
-\frac{\sigma}{2}(g+\frac{1}{6}f^\lambda _{;\lambda })\nu_{\beta }a^{\alpha \beta}\var{u}_\alpha 
-\frac{\sigma}{10}(f+\frac{1}{12}g^\lambda _{;\lambda })\nu_{\beta }a^{\alpha \beta}\var{\psi }_\alpha +\frac{1}{12}g^\alpha \nu_\alpha \var{u}\Bigr] \dd{s}.
\end{multline}
Here, square brackets enclosing indices denote anti-symmetrization: $t^{[\alpha \beta ]}=\frac{1}{2}(t^{\alpha \beta }-t^{\beta \alpha })$. Vector $\nu_\alpha$ represents the outward surface normal to $\partial \mathcal{S}$, and $\mathrm{d}s$ is the length element.  We assume the remaining portion of the shell edge, denoted by $\partial_s$, is free. Due to the arbitrariness of the variations $\var{u}_\alpha$, $\var{\psi}_\alpha$, and $\var{u}$ in $\mathcal{S}$ and on $\partial_s$, we obtain the following second-order partial differential equations from Eq.~\eqref{var3}:
\begin{align}
-t^{\alpha \beta }_{;\beta }+b^\alpha_\beta q^\beta &=f^\alpha -Hg^\alpha +\frac{1}{2}g^\beta b^\alpha_\beta +\frac{\sigma}{2}(g+\frac{1}{6}f^\lambda _{;\lambda })_{;\beta }a^{\alpha \beta}+\frac{1}{12}g^\beta b^\alpha_\beta , \label{equil}
\\
m^{\alpha \beta}_{;\beta }+q^\alpha &=\frac{\sigma}{10}(f+\frac{1}{12}g^\lambda _{;\lambda })_{;\beta }a^{\alpha \beta}+\frac{1}{12}g^\alpha , \label{equilm}
\\
-q^\alpha _{;\alpha }-b_{\alpha \beta} n^{\alpha \beta } &=f-Hg-\frac{1}{2}g^\alpha _{;\alpha }+\sigma H(g+\frac{1}{6}f^\lambda _{;\lambda }), \label{equilq}
\end{align}
where $t^{\alpha \beta }=n^{\alpha \beta }+m^{\lambda [\alpha }b^{\beta ]}_\lambda $. These equations are subject to the kinematic boundary conditions (Eq.~\eqref{kinbc}) on $\partial_k$ and the following natural boundary conditions on $\partial_s$:
\begin{align}
t^{\alpha \beta }\nu _{\beta }&=-\frac{\sigma}{2}(g+\frac{1}{6}f^\lambda _{;\lambda })\nu_{\beta }a^{\alpha \beta}, \label{bc}
\\
-m^{\alpha \beta }\nu_{\beta }&=-\frac{\sigma}{10}(f+\frac{1}{12}g^\lambda _{;\lambda })\nu_{\beta }a^{\alpha \beta}, \label{bcm}
\\
q^\alpha \nu_{\alpha }&=\frac{1}{12}g^\alpha \nu_\alpha . \label{bcq}
\end{align}
The derivation of boundary conditions for a simply supported edge $\partial_{ss}$ is similar. Equations~\eqref{equil}--\eqref{equilq}, \eqref{kinbc}, and \eqref{bc}--\eqref{bcq} constitute the strong formulation of the problem.

The weak and strong formulations of the refined shell theory involve stress resultants, which are related to the following integral characteristics of the 3D stress field:
\begin{equation}
\label{integrals}
T^{\alpha \beta }=\int_{-1/2}^{1/2}\kappa \mu^\alpha_\lambda \sigma^{\beta \lambda}\dd{\xi},\quad M^{\alpha \beta }=\int_{-1/2}^{1/2}\kappa \mu^{\alpha}_\lambda \sigma^{\beta \lambda}\xi\dd{\xi},
\quad
Q^\alpha=\int_{-1/2}^{1/2}\kappa \sigma^{\alpha 3}\dd{\xi}.
\end{equation}
In these integrals, $\xi=x/h$ represents the rescaled transverse coordinate. The terms $\sigma ^{\beta \lambda}$ and $\sigma^{\alpha 3}$ are the contravariant components of the 3D stress tensor expressed in the rescaled shell coordinate system $\{ x^1,x^2,\xi \}$,  while $\kappa =1-2H\xi+K\xi^2$ and $\mu_\alpha ^\beta =\delta_\alpha^\beta-\xi b_\alpha ^\beta$, where $K$ is the Gaussian curvature of the middle surface. Note that $T^{\alpha \beta }$ and $M^{\alpha \beta }$ are non-symmetric. Based on Berdichevsky's work \cite{berdichevsky1979variational,berdichevsky2009variational}, these integral characteristics relate to the stress resultants as follows:
\begin{align}
T^{(\alpha \beta )}&=n^{\alpha \beta }+\frac{\sigma}{2}ga^{\alpha \beta }+\frac{\sigma}{12} f^\lambda_{,\lambda}a^{\alpha \beta }, \label{resultants}
\\
-M^{(\alpha \beta)}&=m^{\alpha \beta }-\frac{\sigma}{10} fa^{\alpha \beta }-\frac{\sigma}{120}g^\lambda_{,\lambda}a^{\alpha \beta }, \label{resultantsm}
\\
Q^\alpha&=q^\alpha-\frac{1}{12}g^\alpha, \label{resultantsq}
\end{align}
where $f,f^\alpha,g,g^\alpha$ are external forces defined in \eqref{extforces}--\eqref{extforcesg}. Note that the right-hand sides of Eqs.~\eqref{resultants}--\eqref{resultantsq} are the partial derivatives of $\Phi -\Theta$, with $\Theta$ being the density of the external work, with respect to $\gamma_{\alpha \beta}$, $\rho_{\alpha \beta}$, and $\varphi_\alpha$, respectively. In this sense they can be regarded as the total membrane forces, bending moments, and transverse shear forces. Note also that all quantities in these equations are rescaled. In Section 5, we will compute the left-hand sides of Eqs.~\eqref{resultants}--\eqref{resultantsq} using the 3D elasticity solution and compare them with the right-hand sides obtained from the 2D refined shell theory solution to verify the asymptotic accuracy of the predicted stress resultants.

\subsection{Discretization}

The weak formulation (Eqs.~\eqref{virtualwork}, \eqref{virtualextension}, \eqref{virtualbending}, and \eqref{virtualshear}) involves geometric correction terms and covariant derivatives of various co- and contravariant vector and tensor components, complicating direct computation of residual forces and stiffness matrices. Therefore, unlike the first-order shear deformation theory for plates presented in \cite{le2024asymptotically}, this refined shell theory employs a different approach. Leveraging the inherent smoothness of the weak formulation, Automatic Differentiation (AD) is used for symbolic representation. Specifically, the \texttt{Trilinos/Sacado} \cite{sacado} package provides the necessary AD capabilities.

The primal variables, $u$, $u_\alpha$, $\psi_\alpha$, and their corresponding variations $\delta u$, $\delta u_\alpha$, $\delta \psi_\alpha$, are assigned degrees of freedom (d.o.f.) with indices $i_u$, $i_{u_\alpha}$, $i_{\psi_\alpha}$, $i_{\delta u}$, $i_{\delta u_\alpha}$, and $i_{\delta \psi_\alpha}$, respectively.  Defining $\var{J} = \delta W^m + \delta W^b + \delta W^s - \delta \mathcal{A}_{\rm{ext}}$ as the first variation of the energy functional, the residual force and stiffness matrix corresponding to primal d.o.f. $\var{a}$ and $b$ ($a, b = u, u_\alpha, \psi_\alpha$) are computed as:
\begin{align} 
r_{\delta a} &= \dfrac{\partial \left( \delta W^m + \delta W^b + \delta W^s - \delta \mathcal{A}_{\rm{ext}} \right)}{\partial \delta a} = \var{J}.dx \left( i_{\delta a} \right) \label{residualf} \\
K_{\delta a, b} &= \dfrac{\partial r_{\delta a}}{\partial b} = \dfrac{ \partial \left( \delta W^m + \delta W^b + \delta W^s - \delta \mathcal{A}_{\rm{ext}} \right)}{\partial \delta a \partial b} =\var{J}.dx \left( i_{\delta a} \right) .dx \left( i_b \right) \label{eq:k_saca}
\end{align}
We note that, to compute the double derivatives in Eq.~\eqref{eq:k_saca}, the symbolic variables $\delta a$ and $b$ must be of type \texttt{Sacado::Fad::DFad<Sacado::Fad::DFad\linebreak <double>>}. The explicit form of the internal force vector required in the residual (Eq.~\eqref{residualf}) is derived in Appendix B. This derivation utilizes the B-matrix method, resulting in the expressions given in Eqs.~\eqref{internalf1}--\eqref{internalf2}.

\subsection{Isogeometric analysis}

As mentioned at the end of Section 3, accurate determination of the true average normal displacement $\check{u}$ and rotation angle $\varphi_\alpha$ necessitates $C^1$ continuity. Isogeometric analysis (IGA) is well-suited for shell geometry discretization for two primary reasons: (i) Higher-order continuity: IGA readily provides the required $C^1$ continuity through the use of non-uniform rational B-splines (NURBS) shape functions, (ii) Geometric flexibility: Complex geometries can be accurately represented by multiple surface patches connected at interfaces, with continuity preserved using techniques like the bending strip method \citep{Kiendl.etal:10}. Furthermore, the order and granularity of the NURBS basis functions can be conveniently increased via $hpk$-refinement.

NURBS shape functions used in IGA read:
\begin{equation}
\mathbf{S} \left( \xi_1, \xi_2 \right) = \sum_{i=1}^m \sum_{j=1}^n \dfrac{N_i^p \left( \xi_1 \right) N_j^q \left( \xi_2 \right)}{\sum_{k=1}^m \sum_{l=1}^n w_{kl} N_k^p \left( \xi_1 \right) N_l^q \left( \xi_2 \right)} \mathbf{P}_{ij}.
\label{eq:nurbs_surface}
\end{equation}
Here, $N_i^p$ and $N_j^q$ are univariate B-spline basis functions of order $p$ and $q$, respectively, computed using the recursive Cox-de-Boor formula:
\begin{align}\label{boor}
N_i^p(\xi) &= \dfrac{\xi-\xi_i}{\xi_{i+p}-\xi_{i}} N^{p-1}_i(\xi) + \dfrac{\xi_{i+p+1}-\xi_i}{\xi_{i+p+1}-\xi_{i+1}} N^{p-1}_{i+1}(\xi) , 
\\ 
N_i^0(\xi)&=\begin{cases}
1 & \xi_i \leq \xi \leq \xi_{i+1} \\
0 & \text{otherwise}
\end{cases}.
\end{align}

$\{ \mathbf{P}_{ij} \}_{0 \leq i \leq m, 0 \leq j \leq n}$ represents the control point grid, and $\{ w_{ij} \}_{0 \leq i \leq m, 0 \leq j \leq n}$ are the corresponding control weights.  The Cox-de-Boor formula (Eq.~\eqref{boor}) requires a global knot vector. To maintain the local nature of finite elements and enable parallel computation during assembly, Bézier decomposition \citep{Borden.etal:11} is employed.

\subsection{Structure of the finite element code}

The proposed 2D-RST element is implemented within the \texttt{PlateAndShellApplication} extension of the modified Kratos Multiphysics framework kernel \cite{prisma,Dadvand.etal:10}.  The \texttt{IsogeometricPlateAndShellApplication} extension enables analysis using IGA. NURBS multipatch structures and Bézier elements are supported in the \texttt{IsogeometricApplication}, which is publicly available at \cite{igaapp}. For an illustration of the relation between software components, the reader is referred to \cite{le2024asymptotically}.

\section{Numerical examples}
\subsection{Semi-cylindrical shell under internal pressure with freely sliding side edges}

\begin{figure}[htb!]
    \centering
    \includegraphics[width=12cm]{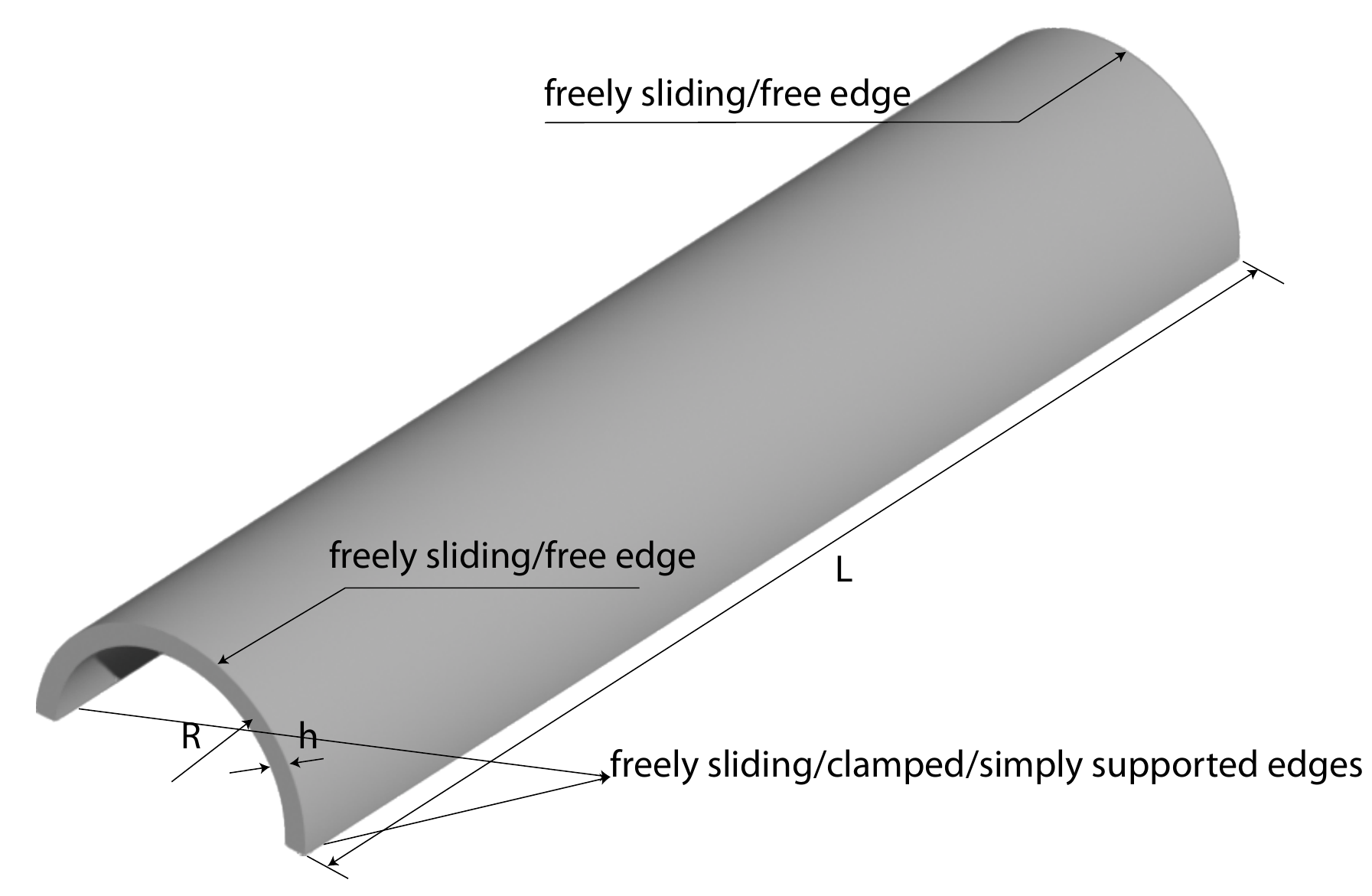}
    \caption{Schematic diagram of the semi-cylindrical shell}
    \label{fig:2}
\end{figure}

We analyze a semi-cylindrical shell occupying the region defined by
\begin{equation}
\label{cylindr}
0 \le x \le L,\quad 0\le \theta \le \pi, \quad R-h/2\le r\le R+h/2
\end{equation}
in cylindrical coordinates $\{ x, \theta ,r\}$ (see Fig.~\ref{fig:2}). The shell is subjected to internal pressure $p$ applied to its inner surface at $r=R-h/2$. For the 2D shell theory, we introduce rescaled surface coordinates $\bar{x}^1=x/h$ varying from $0$ to $\bar{L}=L/h$ and $\bar{x}^2=R \theta/h$ varying from $0$ to $\bar{W}=\pi \bar{R}=\pi R/h$, where $W=\pi R$ is the semi-circular circumference. For brevity, we drop the overbars on rescaled coordinates and quantities in the following theoretical development.

The middle surface in rescaled coordinates is described by:
\begin{equation}
\label{coord}
\vb{z}=-x^1 \vb{e}_1+R \cos \frac{x^2}{R} \vb{e}_2 + R \sin \frac{x^2}{R} \vb{e}_3.
\end{equation}
The tangent vectors to the coordinate lines on the middle surface are:
\begin{equation}
\label{tangent}
\vb{t}_1=\vb{z}_{,1}=-\vb{e}_1,\quad \vb{t}_2=\vb{z}_{,2}=-\sin \frac{x^2}{R} \vb{e}_2 + \cos \frac{x^2}{R} \vb{e}_3.
\end{equation}
The unit normal vector is:
\begin{equation}
\label{normal}
\vb{n}=\frac{\vb{t}_1\cp \vb{t}_2}{\abs {\vb{t}_1\cp \vb{t}_2}}=\cos \frac{x^2}{R} \vb{e}_2 + \sin \frac{x^2}{R} \vb{e}_3.
\end{equation}
In the 2D curvilinear coordinate system $\{ x^1,x^2\}$, the components of the 2D metric tensor are $\delta_{\alpha \beta}$, the Kronecker delta. Consequently, the Christoffel symbols vanish, and covariant derivatives coincide with partial derivatives. Raising or lowering indices does not change tensor and vector components, so we represent them in covariant form with lower indices. The non-zero component of the second fundamental form is:
\begin{equation}
\label{secondq}
b_{11}=0,\quad b_{12}=0,\quad b_{22}=-\frac{1}{R}.
\end{equation}
The mean curvature and $b^\prime _{\alpha \beta}$ are:
\begin{equation}
\label{meanc}
H=-\frac{1}{2R},\quad b^\prime _{11}=\frac{1}{2R}, \quad b^\prime _{12}=0,\quad b^\prime _{22}=-\frac{1}{2R}.
\end{equation}

Using Eqs.~\eqref{scalinggamma}, \eqref{scalingrho}, and \eqref{scalingvarphi}, the extension and bending measures, and the rotation angles due to shear deformation are:
\begin{align}
&\gamma_{11}=u_{1,1},\quad \gamma_{12}=\gamma_{21}= \frac{1}{2}(u_{1,2}+u_{2,1}), \quad \gamma_{22}=u_{2,2}+\frac{u}{R}, \label{scaling1}
\\
&\rho_{11}=-\psi_{1,1},\quad \rho_{22}= -\psi_{2,2}, \label{scaling1rho}
\\
&\rho_{12}= \rho_{21}=-\frac{1}{2}(\psi_{1,2}+\psi_{2,1})+\frac{1}{4R}(u_{1,2}-u_{2,1}), \label{scaling1rho1}
\\
&\varphi_1=u_{,1}+\psi_1 ,\quad \varphi_2= u_{,2}-\frac{u_2}{R}+\psi_2  . \label{scaling1varphi}
\end{align}
The energy density contributions (with overbars dropped) are:
\begin{align}
\Phi_{\rm{cl}}&=\sigma(u_{1,1}+u_{2,2}+\frac{u}{R})^2+(u_{1,1})^2+(u_{2,2}+\frac{u}{R})^2+\frac{1}{2}(u_{1,2}+u_{2,1})^2 \notag
\\
&+\frac{1}{12}\Bigl[ \sigma(\psi_{1,1}+\psi_{2,2})^2+(\psi_{1,1})^2+(\psi_{2,2})^2 \notag
\\
&+\frac{1}{2}\Bigl( -\psi_{1,2}-\psi_{2,1}+\frac{1}{2R}(u_{1,2}-u_{2,1})\Bigr)^2 \Bigr] ,
\end{align}
\begin{align}
\Phi_{\rm{gc}}&=\frac{1}{6R}\Bigl[ \Bigl( \sigma (\frac{6}{5}\sigma-1) -1\Bigr)\rho_{11}\gamma_{11} +(1+\sigma)(1+\frac{6}{5}\sigma)\rho_{22}\gamma_{22} \notag
\\
&+\sigma (1+\frac{6}{5}\sigma)\rho_{22}\gamma_{11}+\sigma(\frac{1}{5}+\frac{6}{5}\sigma)\rho_{11}\gamma_{22} \Bigr] ,\label{phi1}
\end{align}
\begin{equation}
\Phi_{\rm{sc}}=\frac{5}{12}\Bigl[ (u_{,1}+\psi_1)^2+(u_{,2}-\frac{u_2}{R}+\psi_2)^2\Bigr] .
\end{equation}
With pressure acting on the inner face,
\begin{equation}
\label{traction}
f=p, \quad g=-p, \quad f_\alpha=g_\alpha=0,
\end{equation}
and the external work is:
\begin{equation}
\label{work1}
\mathcal{A}_{\rm{ext}}=\int_\mathcal{S} \Bigl[ p(1-\frac{1}{2R})u+\frac{\sigma}{2}p(u_{1,1}+u_{2,2}+\frac{u}{R})-\frac{\sigma}{10}p(\psi_{1,1}+\psi_{2,2})\Bigr] \dd{\omega}.
\end{equation}

We seek analytical solutions as benchmark solutions to validate our finite element implementation of 2D problems. For this purpose, we assume frictionless sliding of the side edges between rigid planes at $x_1=0$ and $x_1=L$.  This implies the following boundary conditions:
\begin{equation}
\label{sidebc}
u_1= \psi_1=0, \quad u_2,\, u,\, \psi_2 \, \text{-arbitrary at $x_1=0,L$}.
\end{equation}
For the bottom edges of the shell, located at $x_2=0,W$, we consider three boundary condition cases:
\begin{enumerate}
  \item Freely sliding edges: 
  \begin{equation}
\label{bc1}
u_1=u_2=0, \quad u \, \text{-arbitrary},\quad \psi_1=\psi_2=0 \quad \text{at $x_2=0,W$}.
\end{equation}
  \item Clamped edges:
  \begin{equation}
\label{bc2}
u_1=u_2=u=\psi_1=\psi_2=0 \quad \text{at $x_2=0,W$}.
\end{equation}
  \item Simply supported edges:
  \begin{equation}
\label{bc3}
u_1=u_2=u=0, \quad \psi_1, \psi_2 \, \text{-arbitrary at $x_2=0,W$}.
\end{equation}
\end{enumerate}

These boundary conditions lead to a state of plane strain, where:
\begin{equation}
\label{plane}
u_1\equiv \psi_1\equiv 0, \quad u_2, \, u, \, \psi_2 \, \text{ are function of $x_2$ only}.
\end{equation}
Consequently, Eqs.~\eqref{scaling1}--\eqref{scaling1varphi} simplify to:
\begin{align}
&\gamma_{11}=0,\quad \gamma_{12}=0,\quad \gamma_{22}\equiv \gamma=u_{2,2} +\frac{1}{R}u, \label{kinematics}
\\
&\rho_{11}=0,\quad \rho_{12}=0,\quad \rho_{22}\equiv \rho=-\psi_{2,2}, \label{kinematicsrho}
\\
&\varphi_1=0,\quad \varphi_2\equiv \varphi=u_{,2} -\frac{1}{R}u_2+\psi_2, \label{kinematicsvarphi}
\end{align}
where the comma before the index 2 denotes the ordinary derivative with respect to $x_2$. The minimization problem (Eq.~\eqref{eq:1}) reduces to a 1-D problem (arc-like model, 1D-RST): Minimize
\begin{multline}
\label{plfunc}
J=\int_{0}^W \Bigl[ (1+\sigma)(u_{2,2}+\frac{1}{R}u)^2+\frac{1+\sigma}{12}(\psi_{2,2})^2-\frac{1}{6}(1+\sigma)(1+\frac{6}{5}\sigma)\frac{1}{R}(u_{2,2}+\frac{1}{R}u)\psi_{2,2} 
\\
+\frac{5}{12}(u_{,2}-\frac{1}{R}u_2+\psi_2)^2-p(1-\frac{1}{2R})u-\frac{\sigma}{2}p(u_{2,2} +\frac{1}{R}u)+\frac{\sigma}{10}p\psi_{2,2} \Bigr] \dd{x_2}
\end{multline}
with respect to $u_2,u,\psi_2$ subject to the kinematic boundary conditions. The integration over $x_1$ yields a constant factor $L$, which is omitted in Eq.~\eqref{plfunc}.

The Euler-Lagrange equations, obtained from the vanishing first variation of functional \eqref{plfunc}, are:
\begin{align}
-n_{,2} -\frac{1}{R}q&=0, \label{pleq}
\\
m_{,2} +q&=0, \label{pleqm}
\\
-q_{,2} +\frac{1}{R}n&=p(1-\frac{1-\sigma }{2R}), \label{pleqq}
\end{align}
where
\begin{align}
n&=\pdv{\Phi}{\gamma}=2(1+\sigma)\gamma+\frac{1}{6}(1+\sigma)(1+\frac{6}{5}\sigma)\frac{1}{R}\rho, \label{plconst}
\\
m&=\pdv{\Phi}{\rho}=\frac{1}{6}(1+\sigma)\rho +\frac{1}{6}(1+\sigma)(1+\frac{6}{5}\sigma)\frac{1}{R}\gamma, \label{plconstm}
\\
q&=\pdv{\Phi}{\varphi}=\frac{5}{6}\varphi. \label{plconstq}
\end{align}
Substituting Eqs.~\eqref{plconst}--\eqref{plconstq} with $\gamma,\rho,\varphi$ from Eqs.~\eqref{kinematics}--\eqref{kinematicsvarphi} into Eqs.~\eqref{pleq}--\eqref{pleqq} yields three second-order ordinary differential equations for $u_2,\psi_2,u$.

For comparison, we also analyze the shell using Sanders-Koiter classical shell theory (CST). Assuming plane strain, the problem reduces to minimizing the 1-D functional (1D-CST)
\begin{equation}
\label{cstfunc}
J_{\text{cl}}=\int_{0}^W \Bigl[ (1+\sigma)(u_{2,2}+\frac{1}{R}u)^2+\frac{1+\sigma}{12}(u_{,22}-\frac{1}{R}u_{2,2})^2
-pu \Bigr] \dd{x_2}
\end{equation}
with respect to $u_2$ and $u$ subject to the kinematic boundary conditions. The Euler-Lagrange equations are:
\begin{align}
-n_{,2} -\frac{1}{R}m_{,2}&=0, \label{cstpleq}
\\
m_{,22} +\frac{1}{R}n&=p, \label{cstpleqm}
\end{align}
where
\begin{align}
n&=2(1+\sigma)\gamma,\quad \gamma=u_{2,2} +\frac{1}{R}u, \label{cstplconst}
\\
m&=\frac{1}{6}(1+\sigma)\rho ,\quad \rho=u_{,22} -\frac{1}{R}u_{2,2}. \label{cstplconstm}
\end{align}

\begin{figure}[htb]
    \centering
    \includegraphics[scale=0.2]{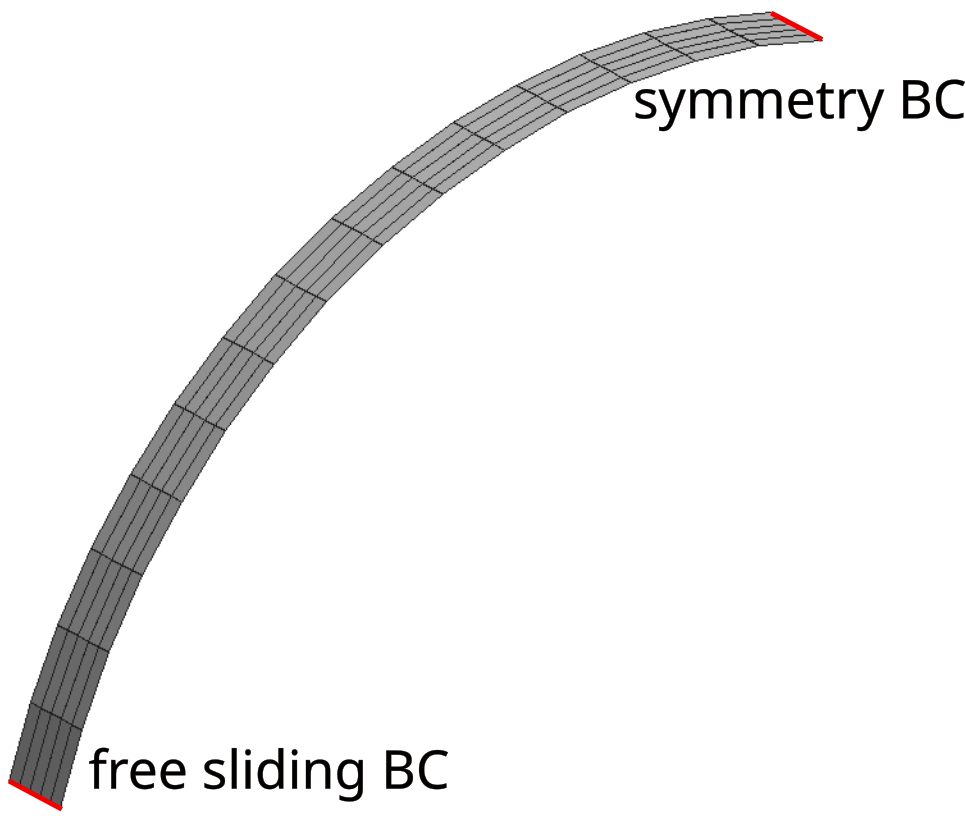}
    \includegraphics[]{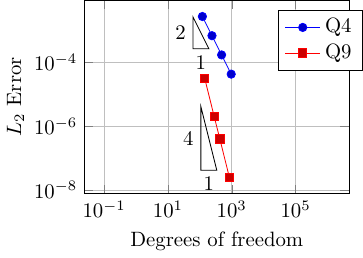}
    \caption{(left) Exemplary mesh for convergence study of the first case, (right) Convergence of 2D-RST displacement in the first case against the analytical solution of 1D-RST in terms of $L_2$ error according to Eq.~\eqref{L2error}.}
    \label{fig:case1}
\end{figure}

For Case 1 (Eq.~\eqref{bc1}), the boundary conditions derived from the vanishing first variation of $J$ and the arbitrariness of $u$ at at $x_2=0,W$ are:
\begin{equation}
\label{1Dbc3}
u_2=\psi_2=0,\quad q=0 \quad \text{at $x_2=0,W$}.
\end{equation}
The solution is
\begin{equation}
\label{sol1}
u=\frac{pR^2}{2(1+\sigma)}(1-\frac{1-\sigma}{2R}),\quad u_2=\psi_2=0,
\end{equation}
yielding
\begin{equation}
\label{exben}
\gamma=\frac{pR}{2(1+\sigma)}(1-\frac{1-\sigma}{2R}),\quad \rho=\varphi=0.
\end{equation}
This solution satisfies Eqs.~\eqref{pleq}, \eqref{plconst}, and the boundary conditions ~\eqref{1Dbc3}. This problem has a known 3D elasticity solution. The radial displacement is \cite{love2013treatise}:
\begin{equation}
\label{radialw}
w_r=\frac{pR^2}{4}(1-\frac{1}{2R})^2\Bigl[ (1-2\nu)(1+\zeta )+(1+\frac{1}{2R})^2\frac{1}{1+\zeta}\Bigr] ,
\end{equation}
where $\zeta=\xi/R$ and $\zeta \in (-1/(2R),1/(2R))$. The average displacement from Eq.~\eqref{radialw} agrees with Eq.~\eqref{sol1} up to $h/R$. To validate the accuracy of the 2D-RST implementation in this Case 1, discretization using linear (Q4) and quadratic (Q9) quadrilateral element are set up for the shell geometry and also to compare the performance. To take advantage of symmetry, the analysis is performed on a quarter-cylinder model, a representative segment of which is depicted in Figure~\ref{fig:case1} (left). The results in Fig.~\ref{fig:case1} (right) shows the convergence of the 2D-RST displacement compared to Eq.~\eqref{sol1} in terms of $L_2$ error defined as follows
\begin{equation}
\label{L2error}
\| e \|_{L_2} = \dfrac{\int_{\mathcal{S}} \| \check{\mathbf{u}} - \check{\mathbf{u}}_{\mathrm{ana}} \|^2 \dd{\omega}}{\int_{\mathcal{S}} \| \check{\mathbf{u}}_{\mathrm{ana}} \|^2 \dd{\omega}}.
\end{equation}
The Q4 element shows a convergence rate of $p+1$ and the Q9 element a rate of $p+2$; the reasons for these particular outcomes are detailed at the end of this subsection.
Classical shell theory (CST) \cite{sanders1959improved,koiter1960consistent}, while accurate for stress resultants and bending moments, exhibits large displacement errors. The 1D-CST solution (Eqs.~\eqref{cstfunc}--\eqref{cstplconstm}) for mean displacements is
\begin{equation}
\label{1DCST}
\begin{split}
&u=-\frac{1}{2}pR^2(1-\nu)\Bigl( \frac{\pi}{2} \sin (x_2/R)-2 \Bigr),
\\
&u_2 = \frac{1}{2}pR^2(1-\nu)\Bigl( \frac{\pi}{2}- \frac{\pi}{2} \cos (x_2/R)-x_2/R \Bigr),
\end{split}
\end{equation}
which results in a $100\%$ error compared to Eqs.~\eqref{sol1} and \eqref{radialw}. This discrepancy highlights the limitations of classical shell theory in accurately predicting displacements, as explained in \cite{berdichevsky1992effect}.

\begin{figure}[htb]
    \centering
    \includegraphics[scale=0.6]{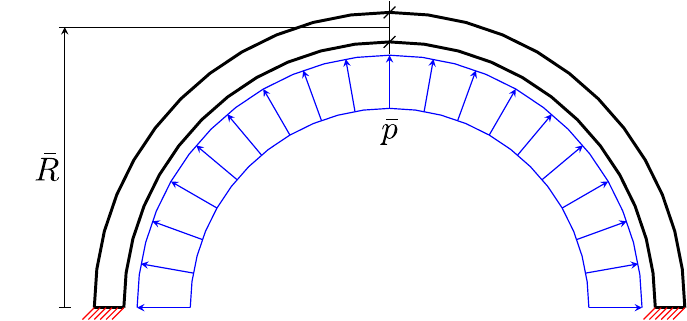}
	\includegraphics[trim=200 180 200 180,clip,scale=0.33]{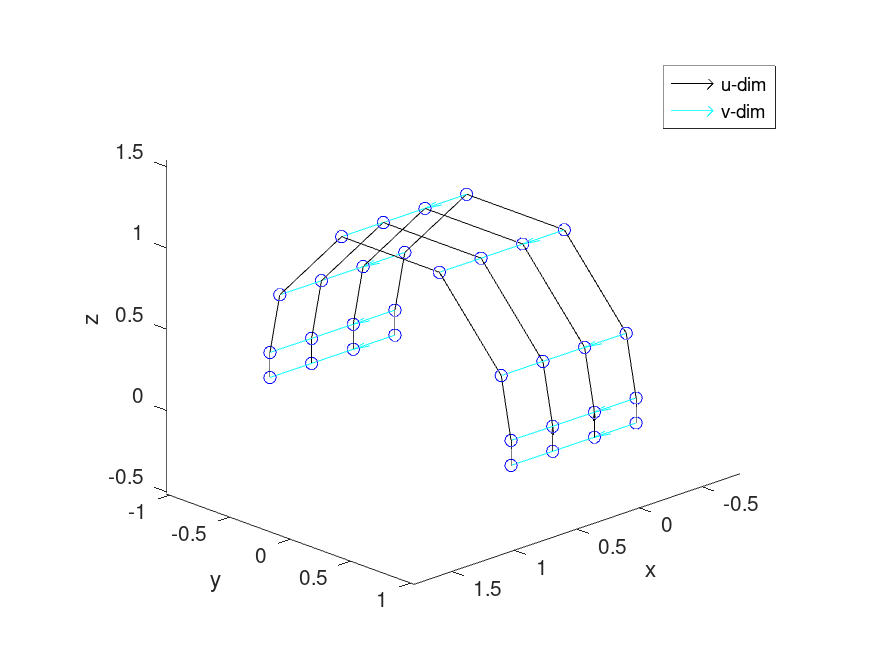}
    \caption{(left) Schematic diagram of the half-ring of thickness $1$ and radius $\bar{R}$ under internal pressure $\bar{p}$, (right) Exemplary control point grid for IGA analysis used in second and third case.}
    \label{fig:halfring}
\end{figure}

For Case 2 with clamped bottom edges, an analytical solution of 1D-RST can be found in terms of exponential functions involving the roots of a cubic polynomial equation. However, this solution is quite cumbersome. Therefore, we opt for a numerical approach. We first rewrite Eqs.~\eqref{kinematics}--\eqref{kinematicsvarphi} as:
\begin{align}
u_{2,2} &=\gamma-\frac{1}{R}u, \label{kinematics1}
\\
\psi_{2,2} &=-\rho, \label{kinematics1rho}
\\
u_{,2} &= \varphi+\frac{1}{R}u_2-\psi_2, \label{kinematics1varphi}
\end{align} 

\begin{table}[htbp] 
    \centering 
    \begin{tabular}{|c|c|c|c|} 
        \hline
        $\bar{L}$ & $\bar{R}$ & $\nu$ & $\bar{p}$ \\
        \hline
        10 & 3 or 10 & 0.3 & 1 \\
        \hline
    \end{tabular}
    \caption{Geometric and material parameters.}
    \label{tab:Parameters} 
\end{table}

\begin{figure}[htb!]
    \centering
    \includegraphics[scale=1.0]{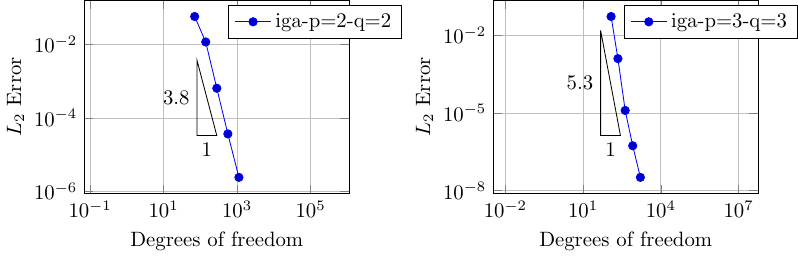}
    \caption{Convergence of 2D-RST displacement in the second case against the benchmark solution of 1D-RST in terms of $L_2$ error \eqref{L2error} for $\bar{R}=10$.}
    \label{fig:case2_u10_conv}
\end{figure}

Next, we eliminate $\rho$ from Eqs.~\eqref{plconst} and \eqref{plconstm} to express $\gamma$ in terms of $n$ and $m$. Differentiating and using Eqs.~\eqref{pleq} and \eqref{pleqm}, we obtain:
\begin{equation}
\label{const1}
\gamma_{,2} =\frac{\sigma }{2(1+\sigma)R[1-\frac{1}{12}(1+\frac{6}{5}\sigma)^2\frac{1}{R^2}]}\varphi.
\end{equation}
Similarly, eliminating $\gamma$ from Eqs.~\eqref{plconst} and \eqref{plconstm}, expressing $\rho$ in terms of $n$ and $m$, differentiating, and using Eqs.~\eqref{pleq} and \eqref{pleqm}, we get:
\begin{equation}
\label{const2}
\rho_{,2} =-\frac{5(1-\frac{1}{12}(1+\frac{6}{5}\sigma)\frac{1}{R^2})}{(1+\sigma)[1-\frac{1}{12}(1+\frac{6}{5}\sigma)^2\frac{1}{R^2}]}\varphi .
\end{equation}
Finally, from Eqs.~\eqref{pleqq} and \eqref{plconstq}, we have:
\begin{equation}
\label{const3}
\varphi_{,2} =\frac{6}{5}[2(1+\sigma)\frac{1}{R}\gamma+\frac{1}{6}(1+\sigma)(1+\frac{6}{5}\sigma)\frac{1}{R^2}\rho-p(1-\frac{1-\sigma}{2R})].
\end{equation}

\begin{figure}[htb!]
    \centering
    \includegraphics[scale=0.86]{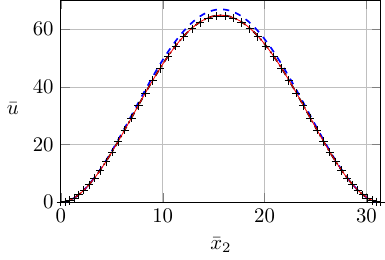}
	\includegraphics[scale=0.86]{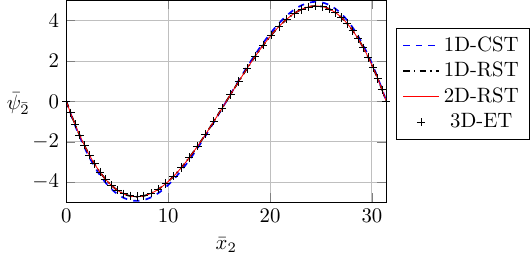}
    \caption{Converged normal displacement $\bar{u}$ (left) and rescaled rotation angle $\bar{\psi}_{\bar{2}}$ (right) of a semi-cylindrical shell under internal pressure with freely sliding side edges and clamped bottom edges (second case) as function of the rescaled circumferential coordinate $\bar{x}_2=\pi \bar{R} \theta$ ($\bar{R}=10$).}
    \label{fig:case2_u10}
\end{figure}

Introducing $v = u_{,2}$ and $\vartheta = \rho_{,2}$, the 1D-CST system (Eqs.~\eqref{cstpleq}--\eqref{cstpleqm} and \eqref{cstplconst}--\eqref{cstplconstm}) becomes:
\begin{align}
u_{,2}&=v, \label{cstodes}
\\
u_{2,2} &=\gamma-\frac{1}{R}u, \label{cstodesu}
\\
\gamma_{,2} &= \frac{1}{12R}\vartheta , \label{cstodesgamma}
\\
v_{,2} &= \rho +\frac{1}{R}(\gamma -\frac{1}{R}u), \label{cstodesv}
\\
\rho_{,2} &= \vartheta , \label{cstodesrho}
\\
\vartheta_{,2} &= \frac{6}{1+\sigma}p-\frac{12}{R}\gamma . \label{cstodesvartheta}
\end{align}

Equations~\eqref{kinematics1}--\eqref{kinematics1varphi} and \eqref{const1}--\eqref{const3} form a system of six first-order ODEs for $u_2, \psi_2, u, \gamma, \rho, \varphi$, subject to the boundary conditions:
\begin{equation}
\label{bc4}
u_2=\psi_2=u=0 \quad \text{at $x_2=0,W$}.
\end{equation}
For 1D-CST, we have six first-order ODEs (Eqs.~\eqref{cstodes}--\eqref{cstodesvartheta}) for $u, u_2, \gamma, v, \rho$, $\vartheta$, subject to:
\begin{equation}
\label{cstbc4}
u=u_2=v=0 \quad \text{at $x_2=0,W$}.
\end{equation}
These two-point boundary-value problems are solved using Matlab's \texttt{bvp4c} function. 

\begin{figure}[htb!]
    \centering
    \includegraphics[scale=0.82]{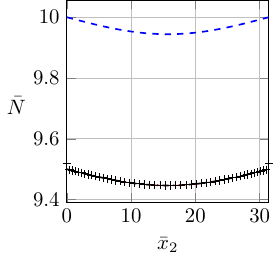}
	\includegraphics[scale=0.82]{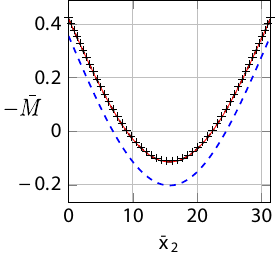}
	\includegraphics[scale=0.82]{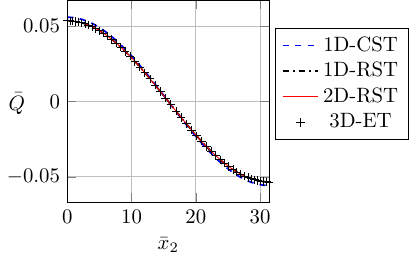}
    \caption{Converged membrane force $\bar{N}$ (left), bending moment $-\bar{M}$ (middle), and shear force $\bar{Q}$ (right) of a semi-cylindrical shell under internal pressure with freely sliding side edges and clamped bottom edges (second case) as function of the rescaled circumferential coordinate $\bar{x}_2=\pi \bar{R} \theta$ ($\bar{R}=10$).}
    \label{fig:case2_nmq10}
\end{figure}

To demonstrate the asymptotic accuracy of our FE-implementation of 2D-RST, we also compare its numerical solution with that of 3D elasticity theory (3D-ET). In this case, under plane strain conditions, the rescaled problem reduces to the two-dimensional analysis of a half-ring having thickness $1$ (see Fig.~\ref{fig:halfring}), which can be solved using 2D isogeometrical analysis (IGA) with solid elements. After solving this 2D problem, the mean radial displacement and rotation angle over the shell's thickness are evaluated in accordance with Eq.~\eqref{eq:ex2_displacement_rotation}. Based on the solution of this problem, we also compute the specific integral characteristics $T^{22}$, $M^{22}$, and $Q^2$ (as defined generally in \eqref{integrals}) to compare them with the corresponding stress resultants in accordance with Eqs.~\eqref{resultants}--\eqref{resultantsq}. For the $O(1/R)$-accuracy of the refined shell theory pertinent to this problem, these integrals can be simplified. We will denote these simplified integral characteristics as:
\begin{equation}
\label{integrals1}
N=\int_{-1/2}^{1/2} \sigma_{\theta \theta}(\xi)\dd{\xi},\quad M=\int_{-1/2}^{1/2}\sigma_{\theta \theta}(\xi)\xi\dd{\xi},
\quad
Q=\int_{-1/2}^{1/2}\sigma_{r\theta}(\xi)\dd{\xi},
\end{equation}
where $\sigma_{r\theta}$ and $\sigma_{\theta \theta}$ represent the physical shear and circumferential stress components, respectively, at the layer defined by $\xi$. To derive these simplified forms, we note the following approximations valid for the current problem. The geometric factor $\kappa =1-2H\xi +K\xi^2$, but since $K=0$, $\kappa=1+\xi/R$. Similarly, the shifter component $\mu^2_2=\delta^2_2-b^2_2\xi = 1+\xi/R$. The contravariant stress component $\sigma^{22}$ is related to the physical circumferential stress $\sigma_{\theta \theta}$ by:
\begin{equation}
\label{sigmaphys}
\sigma^{22}=\frac{1}{g_{22}}\sigma_{\theta \theta},
\end{equation}
where $g_{22}$ is the (2,2)-component of the 3D metric tensor in the shell coordinate system $\{ x^1,x^2,\xi \}$. For this problem, $g_{22}\simeq 1-2b_{22}\xi=1+2\xi/R$. Combining these for the integrand of $T^{22}$ (and $M^{22}$), $\kappa \mu^2_2 \sigma^{22}$ becomes:
\begin{equation}
\label{integrand}
\kappa \mu^2_2 \sigma^{22}=\frac{(1+\xi/R)^2}{1+2\xi/R}\sigma_{\theta \theta} \simeq \frac{1+2\xi/R}{1+2\xi/R}\sigma_{\theta \theta} \simeq \sigma_{\theta \theta},
\end{equation}
as terms of $O(1/R^2)$ are neglected within the RST accuracy. Similarly, for the shear term $Q^2$, the contravariant component $\sigma^{23}$ is related to the physical shear stress component $\sigma_{r \theta}$ by:
\begin{equation}
\label{sigmash}
\sigma^{23}=\frac{1}{\sqrt{g_{22}}}\sigma_{r \theta}=\frac{1}{\sqrt{1+2\xi/R}}\sigma_{r \theta}\simeq \frac{1}{1+\xi/R}\sigma_{r \theta}.
\end{equation}
The integrand for $Q^2$, $\kappa \sigma^{23}$, thus becomes:
\begin{equation}
\label{integrand1}
\kappa \sigma^{23}\simeq \frac{1+\xi/R}{1+\xi/R}\sigma_{r \theta}=\sigma_{r \theta}.
\end{equation}
These approximations justify the simplified forms presented in Eqs.~\eqref{integrals1}.

\begin{figure}[htb!]
    \centering
    \includegraphics[scale=1.0]{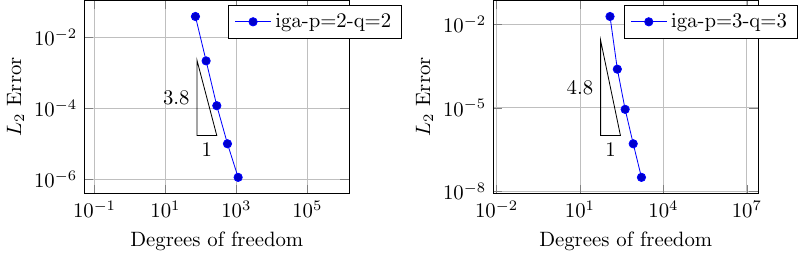}
    \caption{Convergence of 2D-RST displacement in the second case against the benchmark solution of 1D-RST in terms of $L_2$ error \eqref{L2error} for $\bar{R}=3$.}
    \label{fig:case2_u3_conv}
\end{figure}

\begin{figure}[htb]
    \centering
    \includegraphics[scale=0.86]{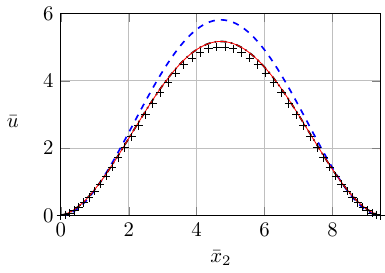}
    \includegraphics[scale=0.86]{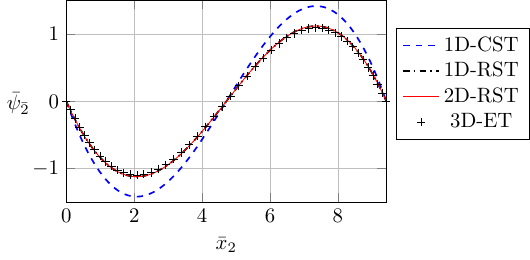}
    \caption{Converged normal displacement $\bar{u}$ and rescaled rotation angle $\bar{\psi}_{\bar{2}}$ of a semi-cylindrical shell under internal pressure with freely sliding side edges and clamped bottom edges (second case) as function of the rescaled circumferential coordinate $\bar{x}_2=\pi \bar{R} \theta$ ($\bar{R}=3$). }
    \label{fig:case2_u3}
\end{figure}

\begin{figure}[htb!]
    \centering
    \includegraphics[scale=0.82]{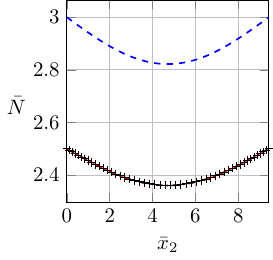}
	\includegraphics[scale=0.82]{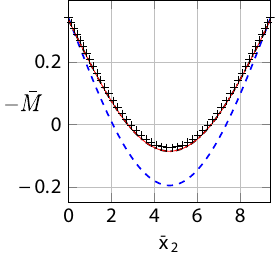}
	\includegraphics[scale=0.82]{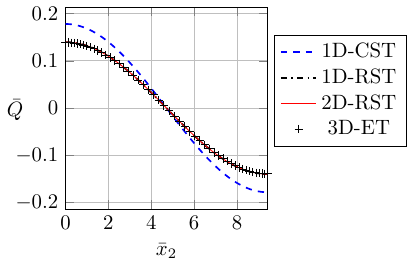}
    \caption{Converged membrane force $\bar{N}$ (left), bending moment $-\bar{M}$ (middle), and shear force $\bar{Q}$ (right) of a semi-cylindrical shell under internal pressure with freely sliding side edges and clamped bottom edges (second case) as function of the rescaled circumferential coordinate $\bar{x}_2=\pi \bar{R} \theta$ ($\bar{R}=3$).}
    \label{fig:case2_nmq3}
\end{figure}

Both the shell analysis and the elastic solid analysis employ cubic-order NURBS discretization with sufficiently fine meshes, by which an exemplary control point grid for the shell is visualized in Fig.~\ref{fig:halfring} (right). The circular shell is modelled using a single NURBS patch. This patch is then applied for both Case 2 and 3. Employing overbars to denote rescaled quantities, we present the results of numerical simulation below. The geometric and material parameters used for the simulations in Sections 5.1 and 5.2 are summarized in Table 1. To normalize the results, the load parameter is set to $\bar{p}=1$. Due to the problem's linearity, the solution for any other load magnitude can then be obtained by simple scaling.

\begin{figure}[htb!]
    \centering
    \includegraphics[scale=1]{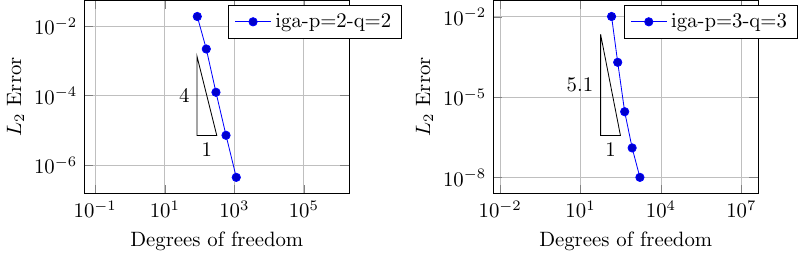}
    \caption{Convergence of 2D-RST displacement in the third case against the benchmark solution of 1D-RST in terms of $L_2$ error \eqref{L2error} for $\bar{R}=10$.}
    \label{fig:case3_u10_conv}
\end{figure}

\begin{figure}[htb!]
    \centering
    \includegraphics[scale=0.86]{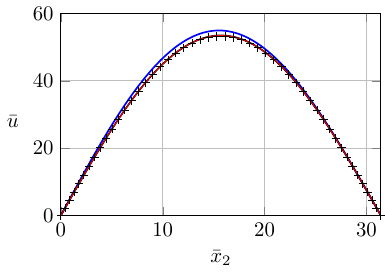}
    \includegraphics[scale=0.86]{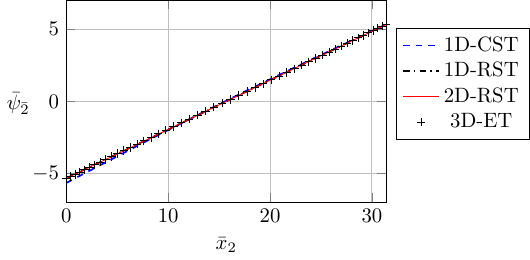}
    \caption{Converged normal displacement $\bar{u}$ (left) and rescaled rotation angle $\bar{\psi}_{\bar{2}}$ (right) of a semi-cylindrical shell under internal pressure with freely sliding side edges and simply-supported bottom edges (third case) as function of the rescaled circumferential coordinate $\bar{x}_2=\pi \bar{R} \theta$ ($\bar{R}=10$).}
    \label{fig:case3_u10}
\end{figure}

\begin{figure}[htb!]
    \centering
    \includegraphics[scale=0.79]{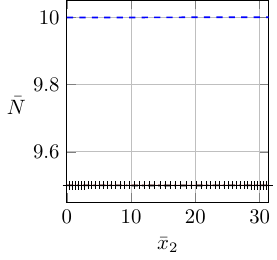}
	\includegraphics[scale=0.79]{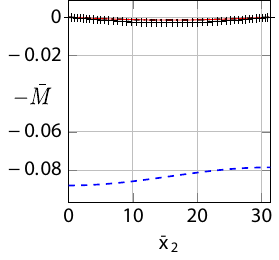}
	\includegraphics[scale=0.79]{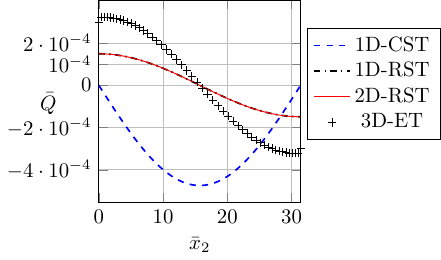}
    \caption{Converged membrane force $\bar{N}$ (left), bending moment $-\bar{M}$ (middle), and shear force $\bar{Q}$ (right) of a semi-cylindrical shell under internal pressure with freely sliding side edges and simply-supported bottom edges (third case) as function of the rescaled circumferential coordinate $\bar{x}_2=\pi \bar{R} \theta$ ($\bar{R}=10$).}
    \label{fig:case3_nmq10}
\end{figure}

Fig.~\ref{fig:case2_u10_conv} shows the convergence rate of the 2D-RST solution to that of the 1D-RST for the thin shell with $\bar{R}=10$, demonstrating that our FE-implementation is free from membrane and shear locking. Fig.~\ref{fig:case2_u10} (left) presents the normal mean displacement $\bar{u}$ versus $\bar{x}_2$ computed using all four theories (1D-CST, 1D-RST, 2D-RST, 3D-ET). The 2D-RST solution (red line) matches that of 1D-RST (black dotted-dashed line) perfectly. Since the shell is thin, this 2D-RST solution agrees quite well with that of 3D-ET (black plus points). The 1D-CST solution deviates noticeably from the others. Fig.~\ref{fig:case2_u10} (right) shows the plot of the rescaled rotation angle $\bar{\psi}_{\bar{2}}$ versus $\bar{x}_2$ computed using all four theories. For classical shell theory, the rotation angle $\bar{\psi}_{\bar{2}}$ is due to bending only and is computed as $\bar{\psi}_{\bar{2}}=-\bar{u}_{,\bar{2}\bar{2}} + \bar{u}_{\bar{2},\bar{2}}/\bar{R}$. We observe a perfect match between the 2D-RST and 1D-RST solutions, and an almost perfect match between them and the 3D-ET solution. Here too, the 1D-CST solution deviates noticeably from the others.

Next, we compare the rescaled stress resultants obtained from various shell theories with the rescaled integral characteristics of the 3D stress state, as defined by \eqref{integrals}. Fig.~\ref{fig:case2_nmq10} shows plots of $\bar{N}$, $-\bar{M}$, and $\bar{Q}$ versus $\bar{x}_2$; these three quantities are computed from the 3D elasticity solution using the simplified integrals presented in \eqref{integrals1}. For comparative purposes, Fig.~\ref{fig:case2_nmq10} also displays: (i) stress resultants from the 2D-RST, as defined by the right-hand sides of Eqs.~\eqref{resultants}--\eqref{resultantsq}; (ii) corresponding quantities from the 1D-RST (specifically, $\bar{n}-\sigma \bar{p}/2$, $\bar{m}-\sigma \bar{p}/10$, $\bar{q}$, where $\bar{n}$, $\bar{m}$, and $\bar{q}$ are from \eqref{plconst}); and
(iii) quantities from the 1D-CST (namely, $\bar{n}$, $\bar{m}$, and $\bar{q}=\bar{m}_{,\bar{2}}$ with $\bar{n}$ and $\bar{m}$ from \eqref{cstplconst} and \eqref{cstplconstm}).
 
The results presented in Fig.~\ref{fig:case2_nmq10} demonstrate that all three stress resultants of the 2D-RST accurately coincide with those of the benchmark 1D-RST. The shear force $\bar{Q}$, which is particularly sensitive to shear locking, exhibits no spurious oscillations--a common issue indicative of shear locking in standard finite element formulations. This result underscores the robustness of our proposed method in this critical aspect. Furthermore, all computed stress resultants show good agreement with the corresponding integral characteristics derived from 3D elasticity theory. This observation confirms the asymptotic accuracy of our formulation for predicting these essential physical quantities. In contrast, the stress resultants from the 1D-CST exhibit significantly larger discrepancies (on the order of $h/R$) when compared to the 3D-ET solution, further highlighting the improved accuracy of the RST.

We now check whether 2D-RST is applicable to moderately thick shells. For this purpose, we set $\bar{R}=3$, while keeping all other parameters unchanged. Fig.~\ref{fig:case2_u3_conv} shows the convergence rate of the 2D-RST solution to that of the 1D-RST. The plots of the mean normal displacement $\bar{u}$ and rescaled rotation angle $\bar{\psi}_{\bar{2}}$ ($-\bar{u}_{,\bar{2}\bar{2}} + \bar{u}_{\bar{2},\bar{2}}/\bar{R}$ for CST) as functions of $\bar{x}_2$ are shown in Fig.~\ref{fig:case2_u3} (left) and Fig.~\ref{fig:case2_u3} (right) for $\nu = 0.3$, $\bar{R} = 3$, and $\bar{p} = 1$. The curves for 1D-RST and 2D-RST, which match each other perfectly, agree well with that of 3D-ET, showing that 2D-RST is applicable also to moderately thick shells. The classical shell theory shows a large deviation from 3D-ET, indicating that 2D-CST can no longer be used for moderately thick shells.

The plots of the rescaled stress resultants versus $\bar{x}_2$ for $\bar{R}=3$ (with other parameters unchanged) are presented in Fig.~\ref{fig:case2_nmq3}. Once again, the curves for the 1D-RST and 2D-RST show a perfect match. Furthermore, they exhibit good agreement with the 3D-ET results, demonstrating that the 2D-RST is also applicable to moderately thick shells. In contrast, the classical shell theory (CST) displays a significant deviation from the 3D-ET, indicating that the CST is not suitable for analyzing moderately thick shells.

\begin{figure}[htb!]
    \centering
    \includegraphics[scale=1]{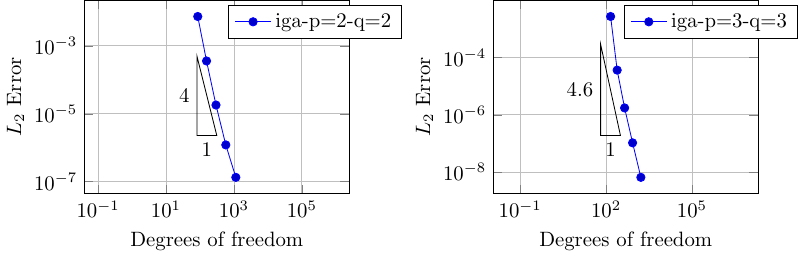}
    \caption{Convergence of 2D-RST displacement in the third case against the benchmark solution of 1D-RST in terms of $L_2$ error \eqref{L2error} for $\bar{R}=3$.}
    \label{fig:case3_u3_conv}
\end{figure}

\begin{figure}[htb!]
    \centering
    \includegraphics[scale=0.86]{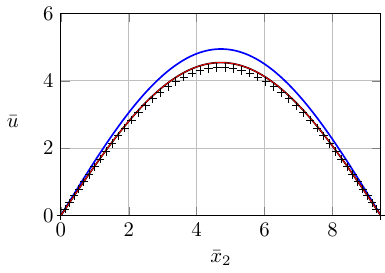}
    \includegraphics[scale=0.86]{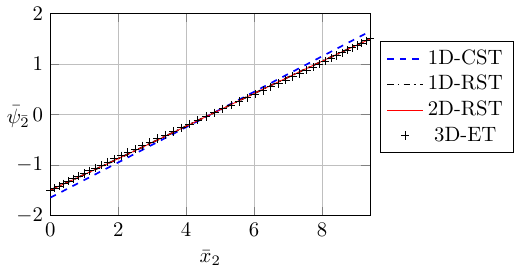}
    \caption{Converged normal displacement $\bar{u}$ and rescaled rotation angle $\bar{\psi}_{\bar{2}}$ of a semi-cylindrical shell under internal pressure with freely sliding side and simply-supported bottom edges (third case) as function of the rescaled circumferential coordinate $\bar{x}_2=\pi \bar{R} \theta$ ($\bar{R}=3$).}
    \label{fig:case3_u3_psi3}
\end{figure}
 
For Case 3 with simply supported bottom edges, the equations remain the same, but the boundary conditions become:
\begin{equation}
\label{bc5}
u_2=u=0, \quad m=\frac{\sigma}{10}p \quad \text{at $x_2=0,W$}.
\end{equation}
The conditions $m = \frac{\sigma}{10}p$ at $x_2 = 0, W$ arise from the vanishing first variation of functional \eqref{plfunc} and the arbitrariness of $\psi_2$ at $x_2 = 0, W$. For 1D-CST, the boundary conditions are:
\begin{equation}
\label{cstbc5}
u=u_2=\vartheta=0 \quad \text{at $x_2=0,W$}.
\end{equation}
These two-point boundary-value problems are solved using Matlab's \texttt{bvp4c} function. Concerning the boundary conditions within 3D-ET: Since exact boundary conditions within 3D-ET for the simply supported edge are absent, we propose boundary conditions that best mimic those of the 2D shell theory. As such, we impose the following constraints: under the plane strain conditions, the mean displacements must vanish:
\begin{equation}
\langle w_2 \rangle =0,\quad \langle w\rangle =0.
\end{equation}
After solving the plane strain problem for the half-ring, the mean radial displacement and rotation angle over the shell's thickness are evaluated in accordance with Eq.~\eqref{eq:ex2_displacement_rotation}, as in the previous case.

\begin{figure}[htb!]
    \centering
    \includegraphics[scale=0.79]{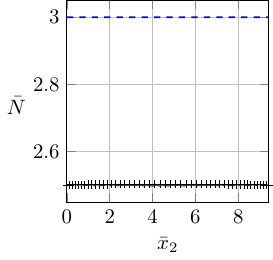}
	\includegraphics[scale=0.79]{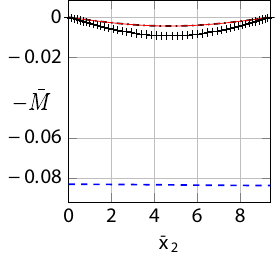}
	\includegraphics[scale=0.79]{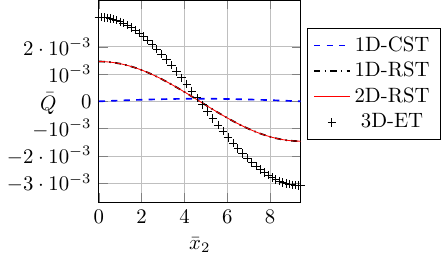}
    \caption{Converged membrane force $\bar{N}$ (left), bending moment $-\bar{M}$ (middle), and shear force $\bar{Q}$ (right) of a semi-cylindrical shell under internal pressure with freely sliding side edges and simply-supported bottom edges (third case) as function of the rescaled circumferential coordinate $\bar{x}_2=\pi \bar{R} \theta$ ($\bar{R}=3$).}
    \label{fig:case3_nmq3}
\end{figure}

The numerical analysis for this third case, featuring simply-supported bottom edges, reinforces the conclusions drawn from the previous scenarios. The 2D-RST demonstrates excellent, locking-free performance and accurately captures the shell's behavior for both thin and moderately thick geometries, while the 1D-CST proves inadequate for moderately thick shells.

For the thin shell ($\bar{R}=10$), the displacement and stress resultant plots (Figs.~\ref{fig:case3_u10} and \ref{fig:case3_nmq10}) once again confirm a perfect match between the 2D-RST and its 1D benchmark, with both showing close agreement with the 3D-ET results. The convergence study in Fig.~\ref{fig:case3_u10_conv} verifies the locking-free nature of the implementation. A notable difference in this case is the slightly larger discrepancy in the shear force when compared to the 3D-ET solution. This is attributed to the approximated boundary conditions used for the 3D solid model, which cannot perfectly replicate a simply-supported edge. As before, the 1D-CST results deviate significantly.

This robust performance is maintained for the moderately thick shell ($\bar{R}=3$), as shown in Figs.~\ref{fig:case3_u3_psi3} and \ref{fig:case3_nmq3}. The RST solutions (1D and 2D) for both displacements and stress resultants align well with the 3D-ET benchmark, confirming the 2D-RST's applicability in this regime. The convergence is shown in Fig.~\ref{fig:case3_u3_conv}. In stark contrast, the 1D-CST results differ substantially, underscoring that classical shell theory is unsuitable for analyzing moderately thick shells.

A key observation from these convergence studies is that the convergence rate in the $L_2$-norm is between $p+1$ and $2p$. We attribute this higher-than-expected rate to a key feature of our asymptotically exact refined shell theory: the inclusion of the bending measure in the determination of the mean normal displacement, $\bar{u}$.

\subsection{Semi-cylindrical shell under internal pressure with free side edges}
 
This subsection investigates the behavior of a semi-cylindrical shell subjected to internal pressure, focusing on the impact of free side edges and varying bottom edge constraints. We analyze the normal displacement $\bar{u}$ and rotation angle $\bar{\psi}_{\bar{2}}$ along the circumferential coordinate, leveraging reduced 1D models and comparisons with 2D-RST solution in the mid-cross-section and plane strain ET solution.

In contrast to the boundary-value problems considered in Section 5.1, these problems do not produce plane strain states. However, for long cylindrical shells with free side edges, where $L \gg R \gg h$, the 1D models and 3D-ET provide accurate benchmark solutions at the shell's mid-cross-section. This is due to the negligible influence of the free side edge boundary effects in the central region, resulting in almost translational invariance along the $x_1$-direction. Consequently, a plane strain state can be assumed there, enabling the application of the arc-like model (1D-RST). Similarly, the 3D-ET simplifies to a 2D problem for the half-ring in the shell's central portion.

\begin{figure}[htb!]
    \centering
	\includegraphics[scale=0.81]{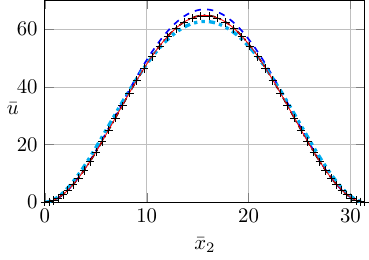}
    \includegraphics[scale=0.81]{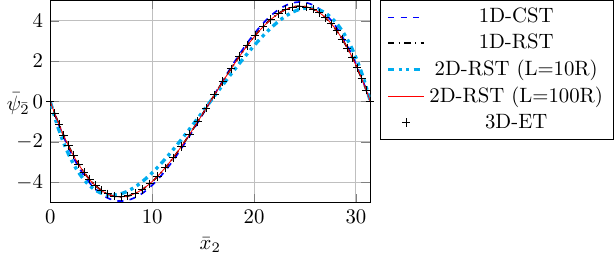}
    \caption{Normal displacement $\bar{u}$ and rescaled rotation angle $\bar{\psi}_{\bar{2}}$ of a semi-cylindrical shell under internal pressure with free side edges and clamped bottom edges (second case) as function of the rescaled circumferential coordinate $\bar{x}_2=\pi \bar{R} \theta$ ($\bar{R}=10$).}
    \label{fig:case2_fe_u10_psi10}
\end{figure}

Fig.~\ref{fig:case2_fe_u10_psi10} illustrates the mean normal displacement $\bar{u}$ and rescaled rotation angle $\bar{\psi}_{\bar{2}}$ as functions of the rescaled circumferential coordinate $\bar{x}_2 = \pi \bar{R} \theta $, calculated using the 2D-RST at the mid-section ($\bar{x}_1 = \bar{L}/2$) for clamped bottom edges. The parameters are taken from Table 1, except for the length $\bar{L}$, which is varied in this analysis. Comparisons are made with 1D-CST, 1D-RST, and plane strain ET solutions. By varying $\bar{L}$ (setting it to $10\bar{R}$ and $100\bar{R}$) while maintaining $\nu = 0.3$, $\bar{R} = 10$, and $p = 1$, we observe that the 2D-RST solution converges to the 1D-RST solution as $\bar{L}$ increases, aligning closely with the 3D-ET results.

\begin{figure}[htb!]
    \centering
	\includegraphics[scale=0.81]{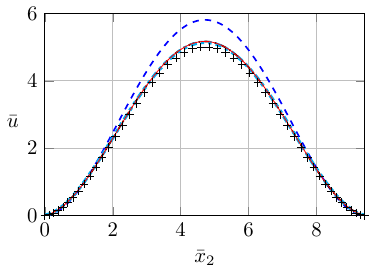}
    \includegraphics[scale=0.81]{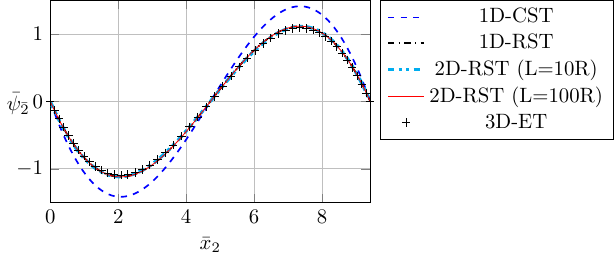}
    \caption{Normal displacement $\bar{u}$ and rescaled rotation angle $\bar{\psi}_{\bar{2}}$ of a semi-cylindrical shell under internal pressure with free side edges and clamped bottom edges (second case) as function of the rescaled circumferential coordinate $\bar{x}_2=\pi \bar{R} \theta$ ($\bar{R}=3$).}
    \label{fig:case2_fe_u3_psi3}
\end{figure}

For moderately thick shells ($\bar{R} = 3$) with clamped bottom edges, Fig.~\ref{fig:case2_fe_u3_psi3} shows similar trends. The 2D-RST solutions at varying $\bar{L}$ values confirm convergence to the 1D-RST solution. Comparisons with 1D-CST and plane strain ET in the mid-cross-section demonstrate the continued applicability of 2D-RST, while highlighting the inadequacy of 2D-CST in this thickness regime.

\begin{figure}[htb!]
    \centering
    \includegraphics[scale=0.79]{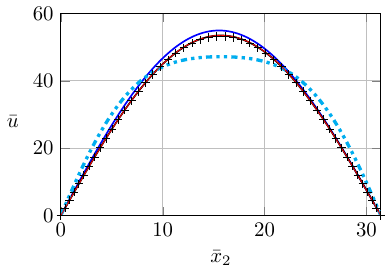}
    \includegraphics[scale=0.79]{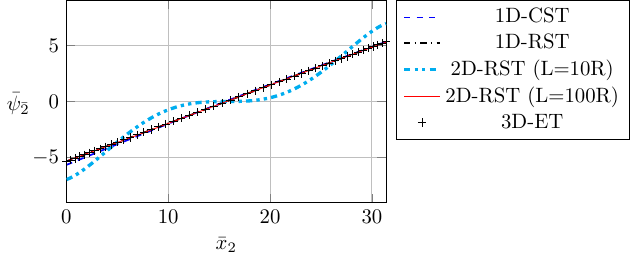}
    \caption{Normal displacement $\bar{u}$ and rescaled rotation angle $\bar{\psi}_{\bar{2}}$ of a semi-cylindrical shell under internal pressure with free side edges and simply-supported bottom edges (third case) as function of the rescaled circumferential coordinate $\bar{x}_2=\pi \bar{R} \theta$ ($\bar{R}=10$).}
    \label{fig:case3_fe_u10_psi10}
\end{figure}

\begin{figure}[htb!]
    \centering
    \includegraphics[scale=0.79]{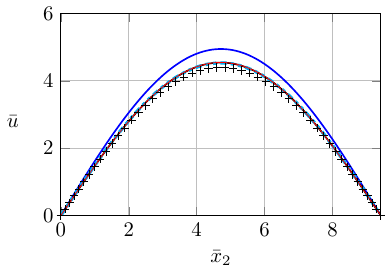}
    \includegraphics[scale=0.79]{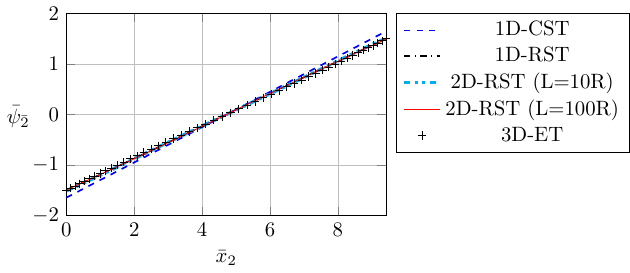}
    \caption{Normal displacement $\bar{u}$ and rescaled rotation angle $\bar{\psi}_{\bar{2}}$ of a semi-cylindrical shell under internal pressure with free side edges and simply-supported bottom edges (third case) as function of the rescaled circumferential coordinate $\bar{x}_2=\pi \bar{R} \theta$ ($\bar{R}=3$).}
    \label{fig:case3_fe_u3_psi3}
\end{figure}

We also examine semi-cylindrical shells with free side edges and simply supported bottom edges. As before, the free side edge effects are negligible in the central region for $\bar{L} \gg \bar{R}$. Numerical results presented in Figs.~\ref{fig:case3_fe_u10_psi10} and \ref{fig:case3_fe_u3_psi3} validate this observation. These figures further illustrate the 2D-RST's applicability to both thin and moderately thick shells, while reinforcing 2D-CST's limitations for moderately thick shells. 

The stress resultant simulations exhibit a similar trend in the mid-section of shells with free side edges; these plots are omitted for brevity. Importantly, the FE implementation of the 2D-RST consistently demonstrates both locking-free behavior and asymptotic accuracy across all numerical simulations.

\subsection{A structure with two plates and a spherical cap under external pressure}
\begin{figure}[htb!]
	\centering
	\includegraphics[scale=0.27,trim=550 110 560 100,clip]{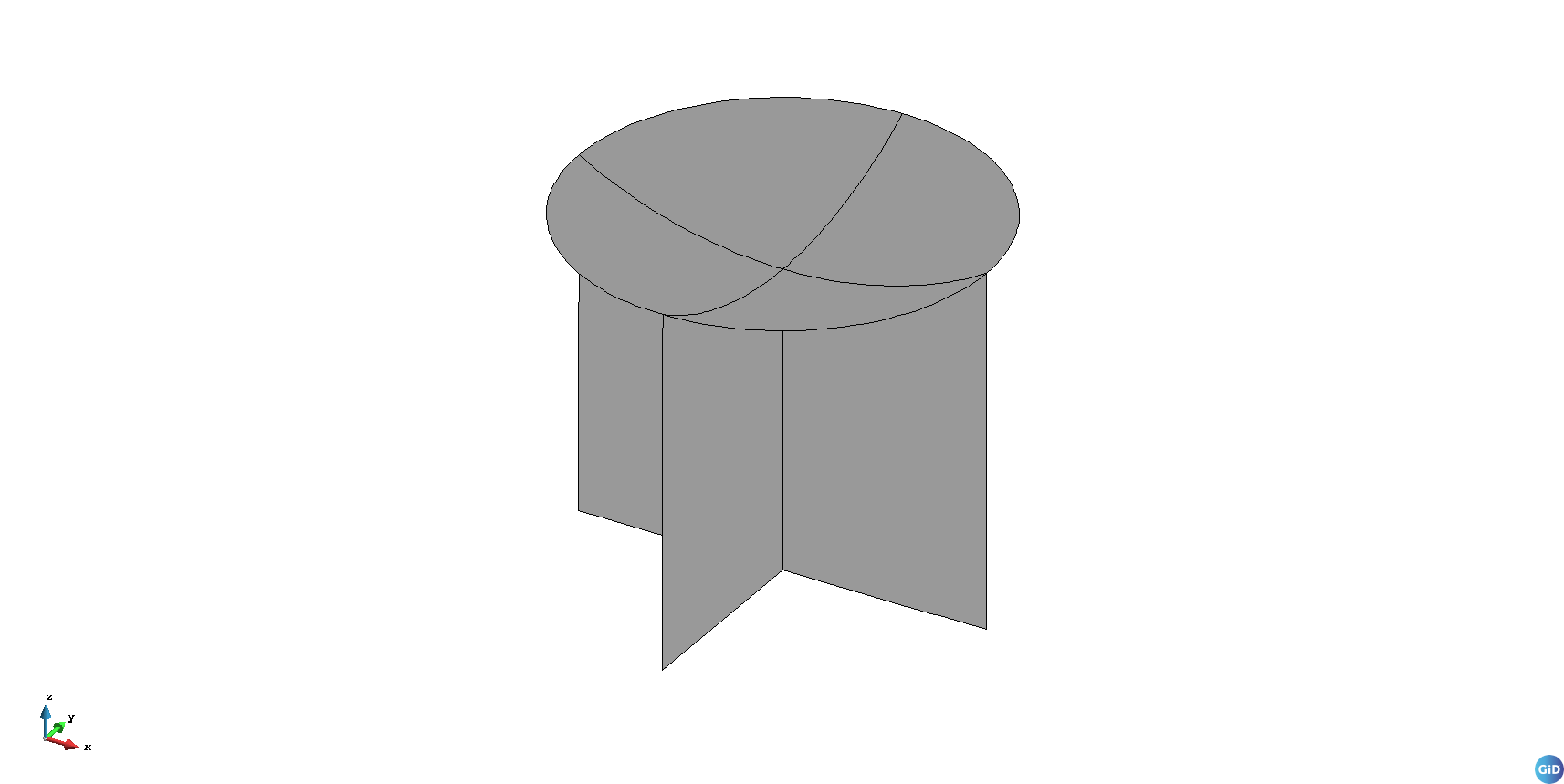}
    \includegraphics[scale=0.35]{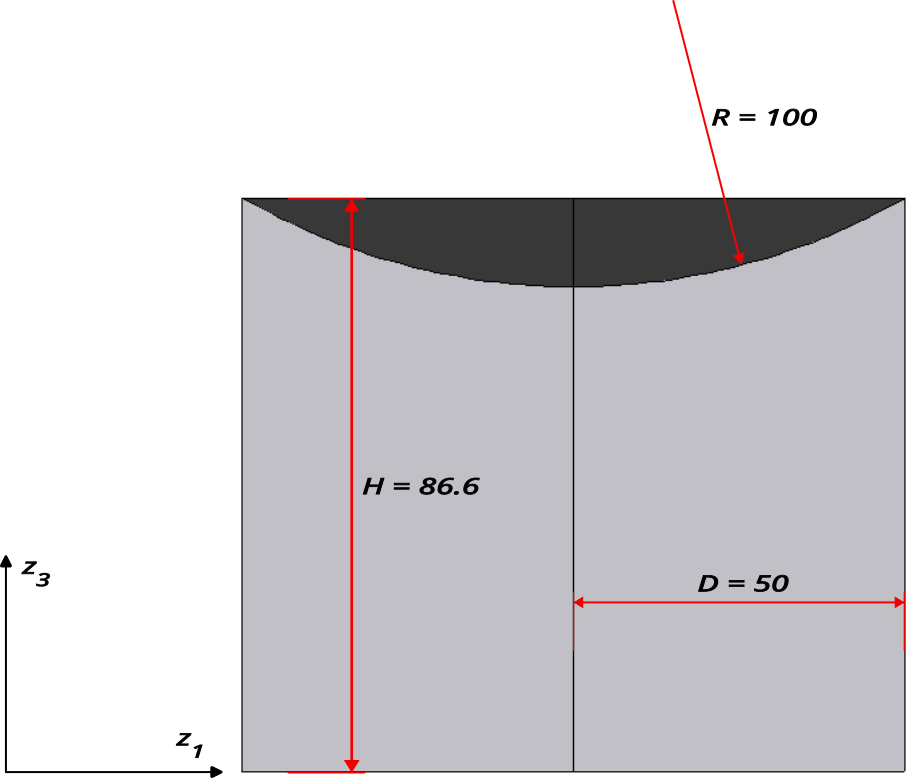}
	\caption{The spherical cap structure and its $(z_1,z_3)$-cross section. The unit of dimension is in cm.}
	\label{fig:structure}
\end{figure}
To showcase the applicability of developed shell formulation in practical application, the deformation of a structure consisting of two orthogonal plates and a spherical cap subjected to external pressure, is investigated. The structure is defined within a Cartesian coordinate system $(z_1, z_2, z_3)$, with the $(z_1, z_2)$-plane coinciding with the base. The mid-surface of the first plate, $\mathcal{P}_1$, is described by:
\begin{equation}
\label{mp1}
\mathcal{P}_1=\{ (z_1,0,z_3)\, | \, \abs{z_1}\le D,\, 0\le z_3\le c-\sqrt{R^2-z_1^2}\},
\end{equation}
where $c = H+\sqrt{R^2-D^2}$ (see Fig.~\ref{fig:structure}). Note that the full width of $\mathcal{P}_1$ extends to $|z_1| \le D + h/2$, but due to $h \ll D$, we focus on the portion intersecting the mid-surface. Similarly, the mid-surface of the second plate, $\mathcal{P}_2$, is:
\begin{equation}
\label{mp2}
\mathcal{P}_2=\{ (0,z_2,z_3)\, | \, \abs{z_2}\le D,\, 0\le z_3\le c-\sqrt{R^2-z_2^2} \} .
\end{equation}
The spherical cap's mid-surface, $\mathcal{P}_3$, is:
\begin{equation}
\label{ms3}
\mathcal{P}_3=\{ (z_1,z_2,z_3)\, | \, z_1^2+z_2^2\le D^2,\, z_3=c-\sqrt{R^2-z_1^2-z_2^2} \}.
\end{equation}
We define $R\phi$ (azimuthal) and $R\theta$ (polar) as the curvilinear coordinates on $\mathcal{P}_3$. All components are made of the same material, with a uniform thickness $h$, where $h \ll D$. The overall mid-surface of the structure is $\mathcal{S} = \mathcal{P}_1 \cup \mathcal{P}_2 \cup \mathcal{P}_3$. A constant pressure $p$ is applied to the cap's top surface, while the plates' bottom edges are clamped. The remaining boundaries are traction-free.

\begin{figure}[htb!]
    \centering
    \includegraphics[scale=0.29,trim=550 110 560 100,clip]{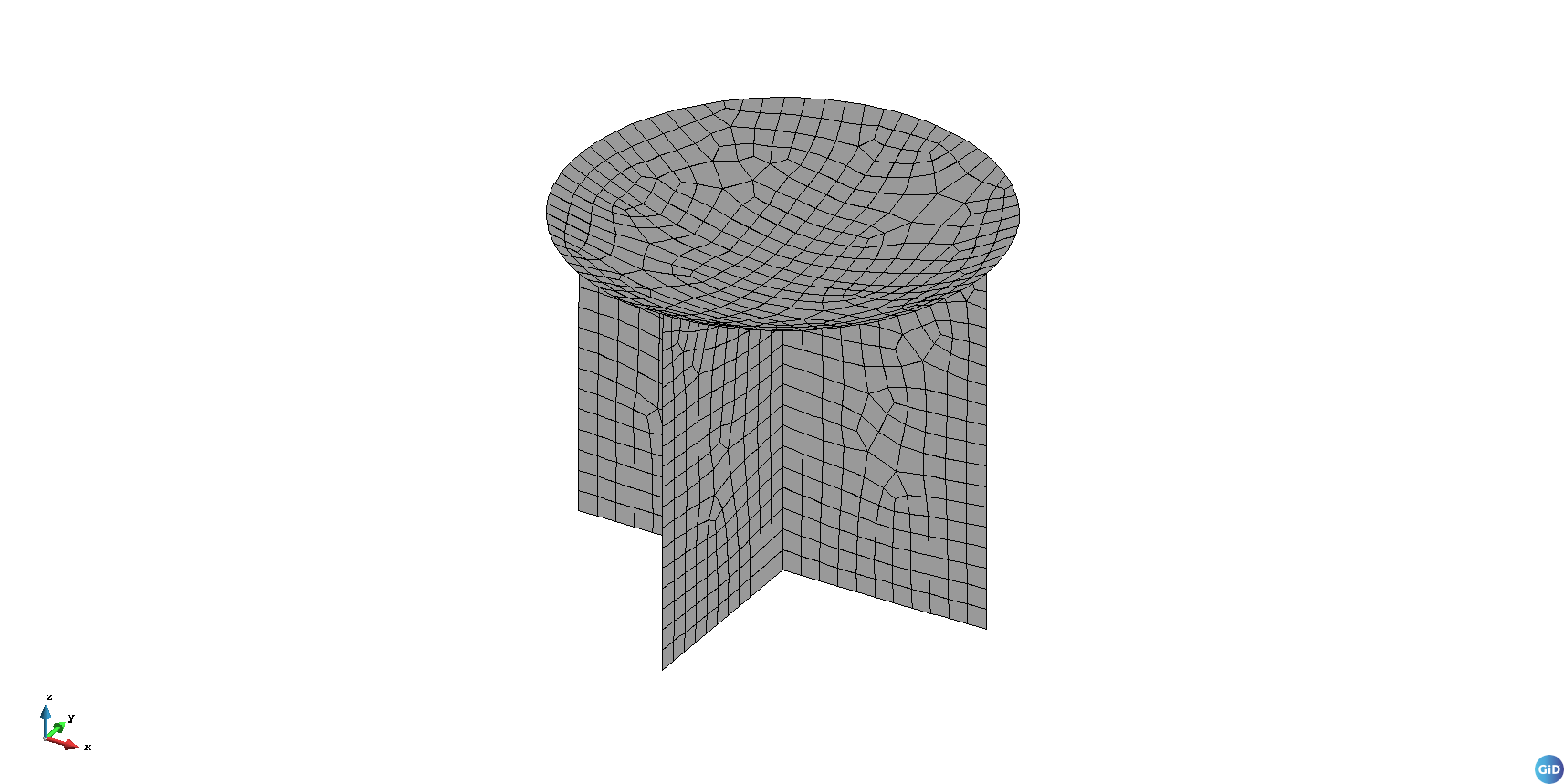}
    \includegraphics[scale=0.29,trim=550 110 560 100,clip]{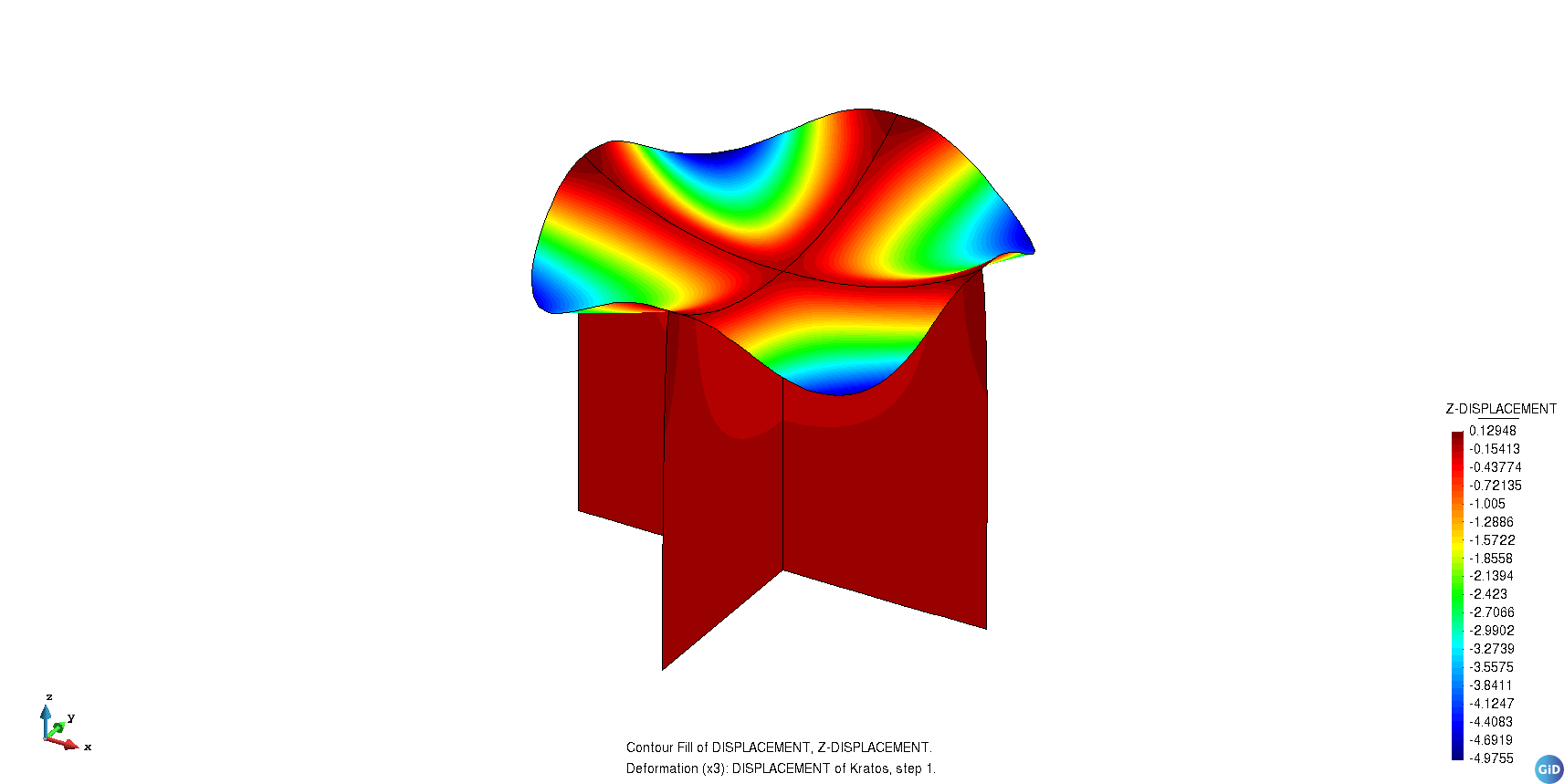}
    \includegraphics[scale=0.36,trim=1500 15 30 445,clip]{disp_z_deformed.png}
    \caption{Mesh (left) and $z_3$-displacement of the spherical cap structure (right). The deformed structure is visualized with scale factor = 3.}
    \label{fig:example}
\end{figure}

The energy functional remains consistent with Eq.~\eqref{energy}, with the energy density $\Phi(\gamma_{\alpha \beta}, \rho_{\alpha \beta}, \varphi_\alpha)$ defined by Eqs.~\eqref{energydensity} and \eqref{cl1}--\eqref{sc1}. The external work is given by:
\begin{equation}\label{work2}
\mathcal{A}_{\rm{ext}}=\int_{\mathcal{P}_3} \Bigl[ -p(1+h/R)u+\frac{\sigma h}{2}p\gamma^\beta_\beta
+\frac{1}{10}\sigma h^2p\rho^\beta_\beta \Bigr] \dd{\omega}.
\end{equation} 
The mid-surfaces $\mathcal{P}_1$, $\mathcal{P}_2$, and $\mathcal{P}_3$ intersect along three lines, where continuity of displacement and rotation vectors is enforced. Rotations about surface normals are constrained, considering only tangential components. At the triple point $(0, 0, c-R)$, where all surfaces meet, the rotation vector is zero.

\begin{figure}[htb!]
    \centering
    \includegraphics[scale=0.27,trim=550 150 550 150,clip]{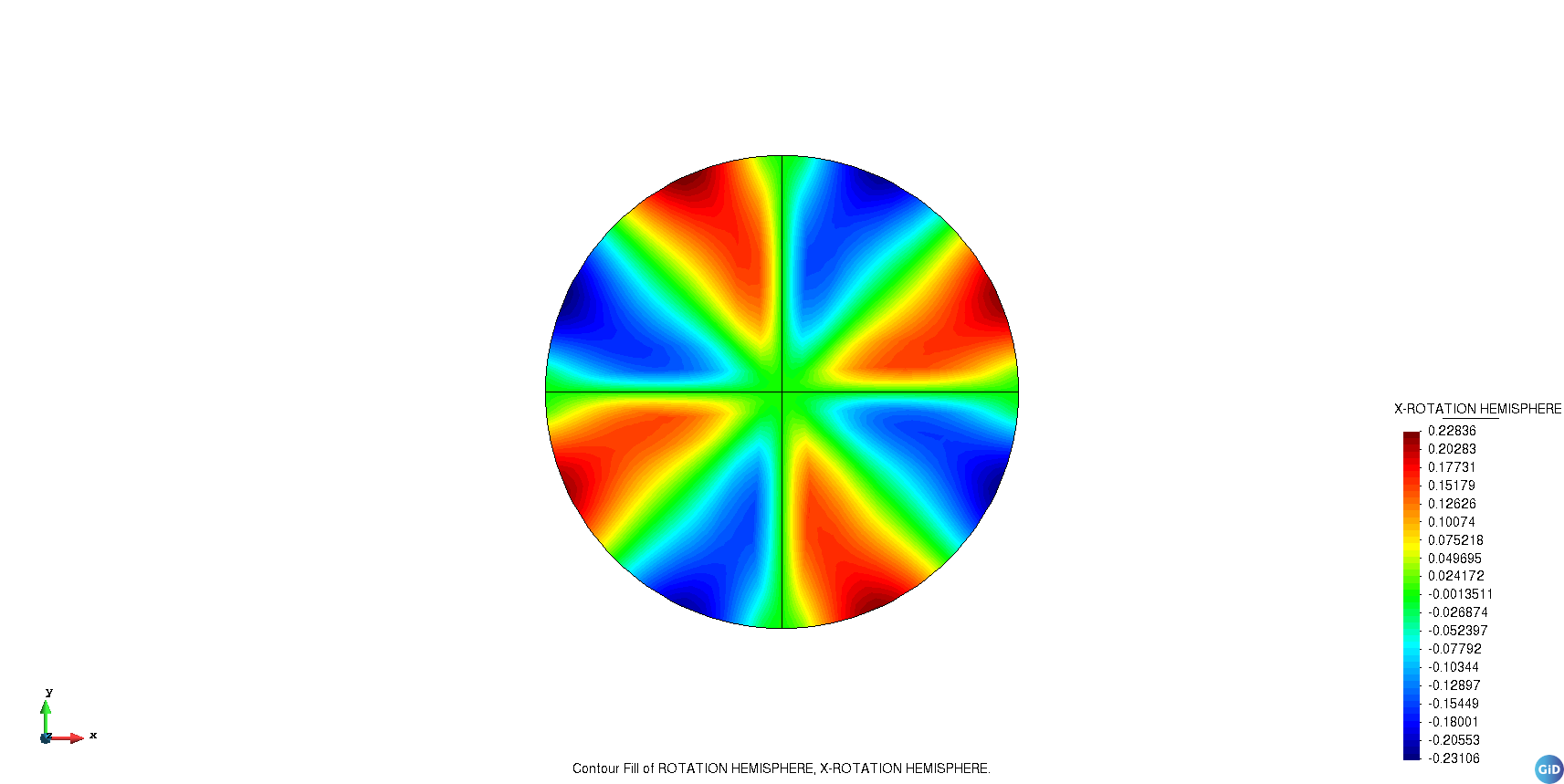}
    \includegraphics[scale=0.33,trim=1450 10 70 445,clip]{rotation_sphere_1.png}
    \includegraphics[scale=0.27,trim=550 150 550 150,clip]{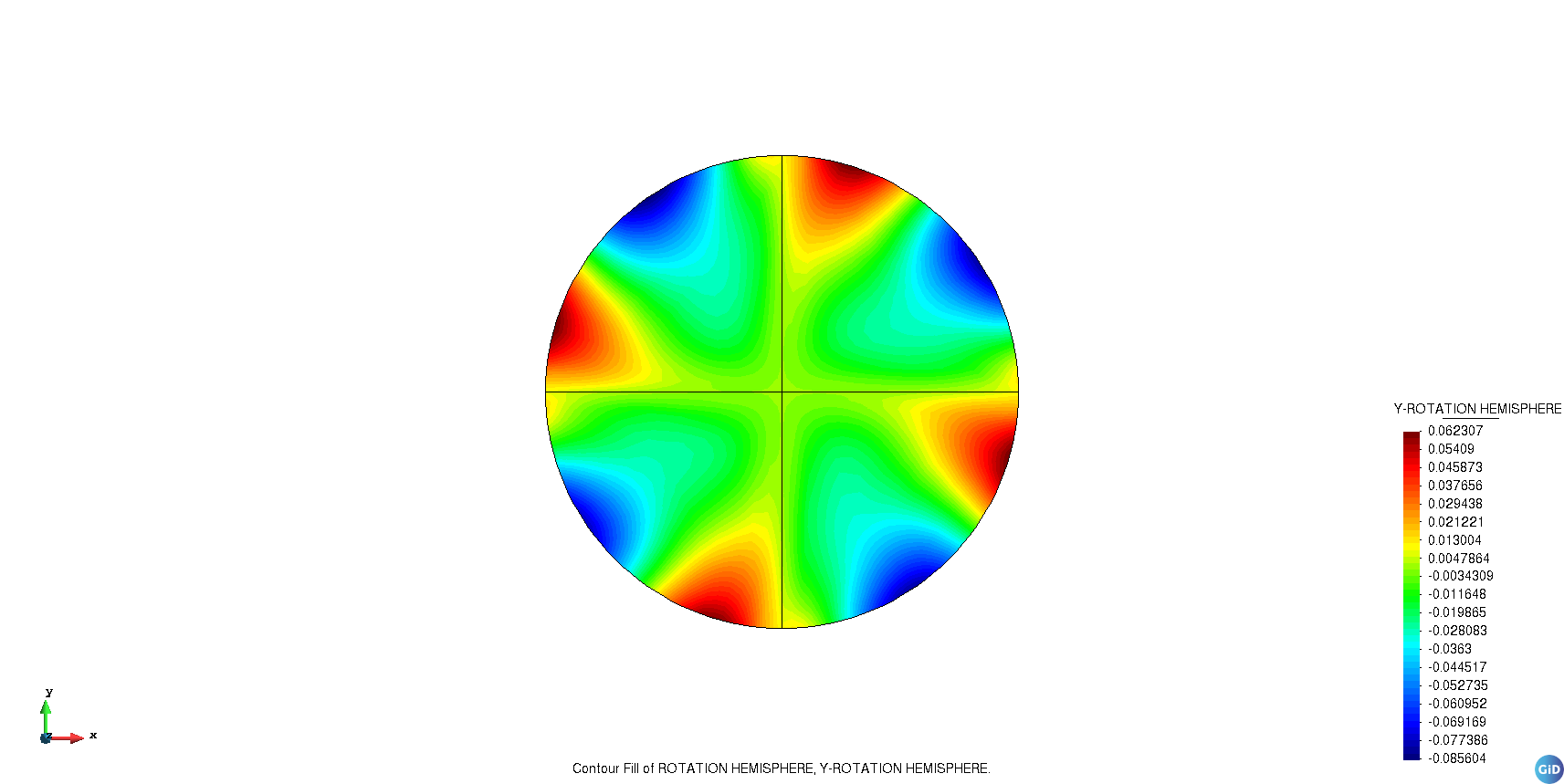}
    \includegraphics[scale=0.33,trim=1450 10 70 445,clip]{rotation_sphere_2.png}
    \caption{Rotation angle projected on $\mathbf{t}_\phi$ (left) and $\mathbf{t}_\theta$ (right) of the spherical cap structure.}
    \label{fig:example_rot}
\end{figure}

The structure is discretized using 1082 Q9 elements with material properties: $E = 2 \times 10^7$ N/cm$^2$, $\nu = 0.3$, and $h = 1$ cm. The applied pressure is $p = 10^3$ N/cm$^2$. Fig.~\ref{fig:example} depicts the deformed structure and its $z_3$-displacement. Fig.~\ref{fig:example_rot} illustrates the rotation vector of the spherical cap projected onto $\mathbf{t}_\phi$ (azimuthal) and  $\mathbf{t}_\theta$ (polar), expressed as functions of the rescaled coordinates $\bar{z}_1$ and $\bar{z}_2$, showing that the azimuthal rotation vanishes along the intersections of the shell and plates.

\section{Conclusion}
This work presents a novel rescaled formulation of a refined shell theory that is inherently free from membrane and shear locking. Combined with an isogeometric finite element implementation, this formulation achieves asymptotic accuracy in shell structure analysis. This combined achievement - asymptotic accuracy and freedom from locking - is key contributions of our paper. Unlike prior methods that often suggest numerical techniques to gain robustness and/or accuracy, hence complicates the formulation further, our approach is mathematical based, concise and rigorous. This enables a straightforward implementation in the numerical code and yields an efficient pathway for analysis of shell structures, especially where accurate stress/displacement prediction in thin/moderately thick shells is crucial. In terms of computational cost, the proposed 2D-RST requires additional resources compared to simpler models. Specifically, its cost is approximately double that of classical theory, a direct result of the added degrees of freedom and the calculation of geometric and shear correction terms. When compared to standard Reissner-Mindlin shell theories, the element-level calculation of geometric correction terms introduces a marginal increase in computational time. This additional expense is not a drawback, but rather a deliberate investment to achieve the significant benefit of asymptotic accuracy. Last but not least, our formulation still relies on the linear kinematics assumption. Future work will develop accurate and locking-free FE implementations for nonlinear refined shell theory, including applications to buckling analysis and shell theories for complex materials.

\section*{Acknowledgements}
\noindent
The second author also received fund­ing from the Marie Curie project TwinSSI (UKRI - Project reference: EP/Z001072/1) and partially from DFG (Project number: 518862444) in collaboration with Vietnam National Foundation for Science and Technology Development (NAFOSTED) under grant number DFG.105-2022.03. These supports are gratefully acknowledged.

\appendix
\section{Derivation of the 2D rescaled variational formulation via an implicit B-matrix}

A large body of literature on the FE-implementation of shells uses the direct B-matrix approach. In this Appendix, we derive the rescaled variational formulation \eqref{eq:1} based on an asymptotically consistent implicit B-matrix to bridge the notational gap with these traditional methods and to clarify the foundation of our proposed formulation. The traditional dimensional reduction, as introduced by Ahmad et al.~\cite{ahmad1970analysis} (see also \cite{hughes2012finite}), begins with a 3D solid shell element and formulates a B-matrix relating strain components to nodal displacements. Kinematic assumptions are then imposed to modify this B-matrix, effectively reducing the 3D elasticity problem to a 2D one. The main drawback of this direct approach is that if the kinematic assumptions are not asymptotically consistent, the FE implementation cannot guarantee asymptotic accuracy. To achieve this, the direct B-matrix approach must be modified.

Our derivation, therefore, begins with the variational principle of 3D elasticity formulated in curvilinear shell coordinates $\{x^1, x^2, x^3 \equiv x\}$~\cite{le1999vibrations}. The objective is to minimize the energy functional
\begin{equation}
\label{3Dfunc}
I[\vb{w}(x^\alpha ,x)]=\int_{\mathcal{S}}\int_{-h/2}^{h/2}W(\vb*{\varepsilon})\kappa \dd{\omega}\dd{x}-\int_{\mathcal{S}_\pm} \tau_i w^i \dd{\omega}_\pm 
\end{equation}
among all kinematically admissible displacement fields $\vb{w}(x^\alpha, x)$. Here, $\kappa =1-2Hx+Kx^2$, $\mathcal{S}_\pm$ are the upper and lower faces of the shell, and $\dd{\omega}_\pm=\kappa_\pm \dd{\omega}$. The strain energy density, $W(\vb*{\varepsilon})$, reads:
\begin{align}
W &=\frac{1}{2}[\lambda (g^{ab}\varepsilon _{ab})^2+
2\mu g^{ac}g^{bd}\varepsilon _{ab}\varepsilon _{cd}] \notag \\
&=\frac{1}{2}[\lambda (g^{\alpha \beta }\varepsilon _{\alpha \beta 
}+\varepsilon _{33})^2 
+2\mu g^{\alpha \gamma }g^{\beta \delta }\varepsilon 
_{\alpha \beta }\varepsilon _{\gamma \delta }+
4\mu g^{\alpha \beta }\varepsilon _{\alpha 3}\varepsilon _
{\beta 3} 
+2\mu \varepsilon ^2_{33}], \label{3Ddensity}
\end{align}
with $g^{ab}$ being the contravariant components of the metric tensor:
\begin{equation}
g^{\alpha \beta }=\frac{1}{\kappa ^2}[(1-2Hx)^2a^{\alpha 
\beta }+2x(1-2Hx)b^{\alpha \beta }+x^2c^{\alpha \beta }],
\,
g^{\alpha 3}=0,\, g^{33}=1,
\label{3Dmetric}
\end{equation}
while $c^{\alpha \beta }=-Ka^{\alpha \beta }+2Hb^{\alpha \beta}$ are the contravariant components of the 2D third quadratic form. The strain components are given in terms of displacements $w_\alpha =t^i_\alpha w_i$, $w=n^iw_i$ by:
\begin{align}
\varepsilon _{\alpha \beta }&=
w_{(\alpha ;\beta )}-b_{\alpha \beta }w-
x b^\lambda _{(\alpha }w_{\lambda ;\beta )}+x c_{\alpha \beta }w,
\label{strain1} \\
2\varepsilon _{\alpha 3}&= 
w_{\alpha ,x}+w_{,\alpha }+b^\lambda _\alpha w_\lambda 
-x b^\lambda _\alpha w_{\lambda ,x},
\label{strain2} 
\\
\varepsilon _{33}&=w_{,x}.
\label{strain3}
\end{align}

The first step is to rescale the coordinates according to
\begin{equation}
\label{scaling3D}
\bar{z}^i=\frac{z^i}{h},\quad \bar{x}^\alpha =\frac{x^\alpha}{h}, \quad \xi =\frac{x}{h},
\end{equation}
which is similar to \eqref{scalingx}. The above variational principle reduces to minimizing the rescaled energy functional:
\begin{equation}
\label{3Dfuncresc}
\bar{I}[\bar{\vb{w}}(\bar{x}^\alpha,\xi)]=\int_{\bar{\mathcal{S}}}\int_{-1/2}^{1/2}\bar{W}(\bar{\vb*{\varepsilon}})\bar{\kappa }\dd{\bar{\omega}}\dd{\xi}-\int_{\bar{\mathcal{S}}_\pm} (\bar{\tau}^{\bar{\alpha}} \bar{w}_{\bar{\alpha}}+\bar{\tau}\bar{w}) \dd{\bar{\omega}}_\pm .
\end{equation}
Here, $\bar{\mathcal{S}}$ is a rescaled domain defined earlier, $\bar{\kappa}=1-2\bar{H}\xi+\bar{K}\xi^2$, $\dd \bar{\omega}$ is the rescaled area element. The rescaled strain energy density is
\begin{equation}
\label{3Denergydensityrs}
\bar{W}=\frac{1}{2}\Bigl[ \frac{2\sigma}{1-\sigma} (\bar{g}^{\bar{\alpha }\bar{\beta }}\bar{\varepsilon }_{\bar{\alpha }\bar{\beta}}+\bar{\varepsilon }_{\bar{3}\bar{3}})^2 
+2 \bar{g}^{\bar{\alpha }\bar{\gamma }}\bar{g}^{\bar{\beta }\bar{\delta }}\bar{\varepsilon} 
_{\bar{\alpha }\bar{\beta }}\bar{\varepsilon }_{\bar{\gamma }\bar{\delta }}+
4 \bar{g}^{\bar{\alpha }\bar{\beta }}\bar{\varepsilon }_{\bar{\alpha }\bar{3}}\bar{\varepsilon }_
{\bar{\beta }\bar{3}} +2 \bar{\varepsilon }^2_{\bar{3}\bar{3}}\Bigr] ,
\end{equation}
where
\begin{equation}
\label{3Dmetricrs}
\bar{g}^{\bar{\alpha }\bar{\beta }}=\frac{1}{\bar{\kappa }^2}[(1-2\bar{H}\xi)^2\bar{a}^{\bar{\alpha }\bar{\beta }}+2\xi(1-2\bar{H}\xi)\bar{b}^{\bar{\alpha }\bar{\beta }}+\xi^2\bar{c}^{\bar{\alpha }\bar{\beta }}].
\end{equation}
Note that the 2D metric tensor does not change, $\bar{a}^{\bar{\alpha }\bar{\beta }}=a^{\alpha \beta }$, while the 2D second and third quadratic forms change according to (see Section 3)
\begin{equation}
\label{bcrs}
\bar{b}^{\bar{\alpha }\bar{\beta }}=hb^{\alpha \beta },\quad \bar{c}^{\bar{\alpha }\bar{\beta }}=h^2c^{\alpha \beta }.
\end{equation}
Thus, $\bar{b}^{\bar{\alpha }\bar{\beta }}$ and $\bar{H}$ are small quantities of order $O(h/R)$, while $\bar{c}^{\bar{\alpha }\bar{\beta }}$ and $\bar{K}$ are small quantities of order $O(h^2/R^2)$. The 3D displacement components remain unchanged,
\begin{equation}
\label{3Ddisplrs}
\bar{w}_{\bar{\alpha }}(\bar{x}^\beta,\xi)=w_\alpha (x^\beta,x),\quad \bar{w}(\bar{x}^\alpha,\xi)=w(x^\alpha,x),
\end{equation}
while the rescaled strain components are:
\begin{align}
\bar{\varepsilon }_{\bar{\alpha }\bar{\beta }}&=
\bar{w}_{(\bar{\alpha };\bar{\beta })}-\bar{b}_{\bar{\alpha }\bar{\beta }}\bar{w}-
\xi \bar{b}^{\bar{\lambda }}_{(\bar{\alpha }}\bar{w}_{\bar{\lambda };\bar{\beta })}+\xi \bar{c}_{\bar{\alpha }\bar{\beta }}\bar{w},
\label{strain1rs} \\
2\bar{\varepsilon }_{\bar{\alpha }\bar{3}}&= 
\bar{w}_{\bar{\alpha },\xi}+\bar{w}_{,\bar{\alpha }}+\bar{b}^{\bar{\lambda }}_{\bar{\alpha }}\bar{w}_{\bar{\lambda}} 
-\xi b^\lambda _\alpha w_{\lambda ,\xi},
\label{strain2rs} 
\\
\bar{\varepsilon }_{\bar{3}\bar{3}}&=\bar{w}_{,\xi}.
\label{strain3rs}
\end{align}
Finally, the rescaled tractions are:
\begin{equation}
\label{tractionrc}
\bar{\tau}_{\bar{\alpha}}=\frac{h}{\mu}\bar{r}^i_{,\bar{\alpha}}\tau_i, \quad \bar{\tau}=\frac{h}{\mu} \bar{n}^i \tau_i .
\end{equation}

Applying the variational-asymptotic analysis to the functional \eqref{3Dfuncresc} (see \cite{berdichevsky2009variational,le1999vibrations}), we derive the following asymptotically consistent Ansatz for the rescaled displacement field:
\begin{multline}
\bar{w}_{\bar{\alpha}} =\bar{u}_{\bar{\alpha}} +\xi \bar{\varphi }_{\bar{\alpha }} -\xi (\bar{u}_{,\bar{\alpha}}+\bar{b}_{\bar{\alpha}}^{\bar{\beta}}\bar{u}_{\bar{\beta}})+\frac{1}{2}\bar{A}^{\bar{\beta}}_{\bar{\beta},\bar{\alpha}}(\xi^2-1/12)-\frac{\sigma}{6}\bar{B}^{\bar{\beta}}_{\bar{\beta},\bar{\alpha}}\xi (\xi^2-\frac{3}{20})
\\
-\frac{5}{3}\tilde{\varphi}_{\bar{\alpha}}\xi (\xi^2-\frac{3}{20})+\frac{5}{6}\bar{g}_{\bar{\alpha}}\xi (\xi^2-3/20)+\frac{1}{2}\bar{f}_{\bar{\alpha}}(\xi^2-1/12),
\label{ldrs} 
\end{multline}
\begin{multline}
\bar{w}=\bar{u}-\xi \sigma \bar{A}^{\bar{\alpha }}_{\bar{\alpha }}
+\frac{1}{2}\sigma \bar{B}^{\bar{\alpha }}_{\bar{\alpha }}(\xi^2-1/12)
\\
+\frac{1-\sigma}{4}\Bigl[ \bar{g} \frac{5}{3} \xi (\xi^2-3/20)+ \bar{f}(\xi^2-1/12) \Bigr] ,
\label{ndrs} 
\end{multline}
with $\bar{f}_{\bar{\alpha}}$, $\bar{g}_{\bar{\alpha}}$, $\bar{f}$, and $\bar{g}$ being defined in \eqref{forcesc}. Functions $\bar{u}_{\bar{\alpha}}$, $\bar{u}$, and $\bar{\varphi }_{\bar{\alpha }}$ are the primary variables of the refined shell theory, representing the mean displacements and rescaled mean rotation angles, which are independent functions of the surface coordinates $\bar{x}^\alpha$. The remaining terms, such as $\bar{A}_{\bar{\alpha } \bar{\beta}}$, $\bar{B}_{\bar{\alpha } \bar{\beta}}$, and $\tilde{\varphi}_{\bar{\alpha}}$, are expressed through these primary variables by the following kinematic relations:
\begin{align}
\bar{A}_{\bar{\alpha } \bar{\beta}}&=\bar{u}_{(\bar{\alpha};\bar{\beta})}-\bar{b}_{\bar{\alpha } \bar{\beta}}\bar{u},\label{Ars}
\\
\bar{B}_{\bar{\alpha } \bar{\beta}}&=\bar{u}_{;\bar{\alpha} \bar{\beta }}+(\bar{u}_{\bar{\lambda}} \bar{b}^{\bar{\lambda }}_{(\bar{\alpha }})_{;\bar{\beta })}
+\bar{b}^{\bar{\lambda }}_{(\bar{\alpha }}\bar{u}_{\bar{\lambda };\bar{\beta })}-\bar{c}_{\bar{\alpha } \bar{\beta}}\bar{u} - \bar{\varphi}_{(\bar{\alpha};\bar{\beta})}, \label{Brs} 
\\
\tilde{\varphi}_{\bar{\alpha}}&=\bar{\varphi }_{\bar{\alpha }}-\frac{\sigma}{60}\bar{B}^{\bar{\beta}}_{\bar{\beta},\bar{\alpha}}.
\end{align}
This displacement Ansatz is equivalent to introducing kinematic constraints in the traditional B-matrix approach; its key advantage is that it is derived assumption-free and is asymptotically consistent.

Based on this Ansatz, the asymptotically accurate rescaled strain components are obtained as:
\begin{align}
\bar{\varepsilon }_{\bar{\alpha }\bar{\beta }}&=
\bar{A}_{\bar{\alpha } \bar{\beta}}-\bar{B}_{\bar{\alpha } \bar{\beta}} \xi+\bar{B}_{\bar{\lambda} (\bar{\alpha }}\bar{b}^{\bar{\lambda}}_{\bar{\beta})}\xi^2,
\label{strain1rs} \\
2\bar{\varepsilon }_{\bar{\alpha }\bar{3}}&= 
\tilde{\varphi}_{\bar{\alpha }}+5\tilde{\varphi}_{\bar{\alpha }}(\xi^2-1/20)+\frac{1}{2}\Bigl[ \bar{g}_{\bar{\alpha }} (5\xi^2-1/4)+ 2\bar{f}_{\bar{\alpha }}\xi \Bigr],
\label{strain2rs} 
\\
\bar{\varepsilon }_{\bar{3}\bar{3}}&=-\sigma \bar{A}^{\bar{\alpha }}_{\bar{\alpha }}+\sigma \bar{B}^{\bar{\alpha }}_{\bar{\alpha }}\xi +\frac{1-\sigma}{4}\Bigl[ \bar{g} (5\xi^2-1/4)+ 2\bar{f}\xi \Bigr].
\label{strain3rs}
\end{align}
The terms $\bar{A}_{\bar{\alpha } \bar{\beta}}$, $\bar{B}_{\bar{\alpha } \bar{\beta}}$, and $\tilde{\varphi}_{\bar{\alpha}}$ contain covariant derivatives along the surface directions. Following the 2D FE-discretization, these derivatives are expressed in terms of the nodal primary variables. Consequently, the system of equations~\eqref{strain1rs}--\eqref{strain3rs} acts as an implicit B-matrix, establishing the required relationship between the strain components and the nodal degrees of freedom.

The final step involves substituting these expressions into the rescaled functional and integrating over the thickness variable $\xi \in (-1/2,1/2)$. Retaining terms up to the appropriate asymptotic order yields the 2D energy functional. Following a procedure similar to that detailed in~\cite{berdichevsky2009variational,berdichevsky1979variational}, this functional is then simplified to the form shown in Eq.~\eqref{eq:1} through a change of variables for the rotations and transverse displacement, and a redefinition of the extensional and bending measures, $\bar{\gamma}_{\bar{\alpha }\bar{\beta }}$ and $\bar{\rho}_{\bar{\alpha }\bar{\beta}}$.

\section{2D rescaled weak formulation in matrix-vector form}

The rescaled variational problem presented in Eq.~\eqref{eq:1} is written in a tensor form that can be directly implemented in finite element code as shown in Section 4. However, presenting the formulation in a matrix-vector form, following the standard B-matrix approach, is advantageous for two primary reasons: it can lead to higher computational efficiency through vectorization, and it facilitates the subsequent linearization required to obtain the stiffness matrix.

This Appendix details the derivation of the explicit B-matrices and the consistent tangent stiffness matrix. For brevity, the overbars on all rescaled quantities are dropped.

Following the weak formulation in Eq.~\eqref{var1}, the virtual work density of the internal forces can be expressed as:
\begin{equation}
    \delta \Phi = n^{\alpha \beta} \delta \gamma_{\alpha \beta} + m^{\alpha \beta} \delta \rho_{\alpha \beta} + q^\alpha \delta \varphi_\alpha .
\end{equation}
For computational purposes, the components of symmetric, rank-2 tensors are often arranged into row vectors using 2D Voigt notation. This procedure is defined differently for stress and strain tensors to ensure their work-conjugate relationship is preserved in vector form. The mapping for stress-like tensors, such as the contravariant stress resultants, $n^{\alpha \beta}$ and $m^{\alpha \beta}$, is defined by the operator $\mathcal{M}$ as:
\begin{equation}
\label{M}
\text{Given } n^{\alpha \beta} = \begin{pmatrix}
n^{11}  & n^{12} \\
n^{12}  & n^{22}
\end{pmatrix},
\quad \text{then} \quad
\mathcal{M}(n^{\alpha \beta}) = (n^{11}, n^{22}, n^{12}).
\end{equation}
For the corresponding strain-like covariant tensors (e.g., $\gamma_{\alpha \beta}$, $\rho_{\alpha \beta}$, and their variations), the mapping introduces a factor of 2 for the shear components:
\begin{equation}
\label{M1}
\mathcal{M}(\gamma_{\alpha \beta}) = (\gamma_{11}, \gamma_{22}, 2\gamma_{12}).
\end{equation}
This convention is essential as it preserves the inner product for virtual work, ensuring that the tensor contraction equals the vector dot product: $n^{\alpha \beta} \delta \gamma_{\alpha \beta}=[\mathcal{M}(n^{\alpha \beta})][\mathcal{M}(\delta \gamma_{\alpha \beta})]^T$. Correspondingly, the variations of the strain measures ($\delta \gamma_{\alpha \beta}$, $\delta \rho_{\alpha \beta}$, and $\delta \varphi_\alpha$) are related to the vector of elemental degrees of freedom, $\delta \mathbf{d}$, through their respective B-matrices:
\begin{equation}
\mathcal{M}(\delta \gamma_{\alpha \beta}) = \delta \mathbf{d}^T \mathbf{B}^n \, , \quad \mathcal{M}(\delta \rho_{\alpha \beta}) = \delta \mathbf{d}^T \mathbf{B}^m, \quad \text{and} \quad \delta \varphi_{\alpha} = \delta \mathbf{d}^T \mathbf{B}^q,
\end{equation}
in which $\mathbf{d}$ is the column vector representing the elemental d.o.f.s
\begin{equation}
\mathbf{d} = \begin{bmatrix} \cdots & \hat{\mathbf{u}}_i & \cdots & \hat{\boldsymbol{\psi}}_j & \cdots \end{bmatrix}^T \, , \quad i=1 \cdots n_u, j=1 \cdots n_\psi .
\end{equation}

Under the linear kinematics assumption, the local coordinates vectors $\mathbf{n}, \mathbf{t}_\alpha$ are assumed to be displacement- and rotation-independent. Therefore, their variations, along with those of other geometric quantities, vanish. With this simplification, the variations of the strain measures become:
\begin{align}
\delta \gamma_{\alpha \beta} &= \delta u_{(\alpha;\beta)} - b_{\alpha \beta} \delta u, \\
\delta \rho_{\alpha \beta} &= -\delta \psi_{(\alpha;\beta)} + b^\lambda_{(\alpha} \delta \varpi_{\beta)\lambda}, \\
\delta \varphi_\alpha &= \delta \psi_\alpha + \mathbf{n} \cdot \delta \mathbf{u}_{,\alpha},
\end{align}
where the constituent variations are defined as:
\begin{align}
&\delta u_{\alpha;\beta} = \delta u_{\alpha,\beta} - \Gamma^\lambda_{\alpha \beta} \delta u_\lambda \, , \quad \delta \psi_{\alpha;\beta} = \delta \psi_{\alpha,\beta} - \Gamma^\lambda_{\alpha \beta} \delta \psi_\lambda , \label{eq:varupsi} \\
&\delta \varpi_{\alpha \beta} = \dfrac{1}{2} \left( \delta u_{\beta,\alpha} - \delta u_{\alpha,\beta} \right), \label{eq:varpi} \\
&\delta \psi_\alpha = \mathbf{t}_\alpha \cdot \delta \boldsymbol{\psi} \, , \quad   \delta u_\alpha = \mathbf{t}_\alpha \cdot \delta \mathbf{u} \, , \quad \delta u = \mathbf{n} \cdot \delta \mathbf{u},
\label{eq:var}
\end{align}
and
\begin{equation}
\delta u_{\alpha,\beta} = \mathbf{t}_{\alpha,\beta} \cdot \delta \mathbf{u} + \mathbf{t}_\alpha \cdot \delta \mathbf{u}_{,\beta} \, , \quad \delta \psi_{\alpha,\beta} = \mathbf{t}_{\alpha,\beta} \cdot \delta \boldsymbol{\psi} + \mathbf{t}_\alpha \cdot \delta \boldsymbol{\psi}_{,\beta}. \label{eq:varder}
\end{equation}

By denoting $\mathbf{N}_u[(dim \times n_u), dim]$ and $\mathbf{N}_\psi[(dim \times n_\psi), dim]$ as the interpolation matrices for the displacement and rotation fields, the variation terms in \eqref{eq:varupsi}-\eqref{eq:varder} can be expressed in terms of the nodal variations $\delta \hat{\mathbf{u}}$ and $\delta \hat{\boldsymbol{\psi}}$:
\begin{align}
\delta u_{\alpha;\beta} &= \delta \hat{\mathbf{u}}^T \left( \mathbf{N}_u \mathbf{t}_{\alpha,\beta} + \mathbf{N}_{u,\beta} \mathbf{t}_{\alpha} - \Gamma^\lambda_{\alpha \beta} \mathbf{N}_u \mathbf{t}_\lambda \right) = \delta \hat{\mathbf{u}}^T \mathbf{m}^u_{\alpha \beta}, \\
\delta \psi_{\alpha;\beta} &= \delta \hat{\boldsymbol{\psi}}^T \left( \mathbf{N}_\psi \mathbf{t}_{\alpha,\beta} + \mathbf{N}_{\psi,\beta} \mathbf{t}_{\alpha} - \Gamma^\lambda_{\alpha \beta} \mathbf{N}_\psi \mathbf{t}_\lambda \right) = \delta \hat{\boldsymbol{\psi}}^T \mathbf{m}^\psi_{\alpha \beta}, \\
\delta \varpi_{\alpha \beta} &= \dfrac{1}{2} \delta \hat{\mathbf{u}}^T \left( \mathbf{N}_u \mathbf{t}_{\beta,\alpha} + \mathbf{N}_{u,\alpha} \mathbf{t}_{\beta} - \mathbf{N}_u \mathbf{t}_{\alpha,\beta} - \mathbf{N}_{u,\beta} \mathbf{t}_{\alpha} \right) = \delta \hat{\mathbf{u}}^T \mathbf{m}^\varpi_{\alpha \beta}, \\
\delta \psi_\alpha &= \delta \hat{\boldsymbol{\psi}}^T \mathbf{N}_\psi \mathbf{t}_\alpha \, , \, \delta u_\alpha = \delta \hat{\mathbf{u}}^T \mathbf{N}_u \mathbf{t}_\alpha \, , \, \delta u = \delta \hat{\mathbf{u}}^T \mathbf{N}_u \mathbf{n}  \, , \, \delta \mathbf{u}_{,\alpha} = \delta \hat{\mathbf{u}}^T \mathbf{N}_{u,\alpha}.
\end{align}
This leads to the final form of the strain measure variations:
\begin{align}
\delta \gamma_{\alpha \beta} &= \delta \hat{\mathbf{u}}^T \left( \mathbf{m}^u_{(\alpha \beta)} - b_{\alpha \beta} \mathbf{N}_u \mathbf{n} \right), \\
\delta \rho_{\alpha \beta} &= - \delta \hat{\boldsymbol{\psi}}^T \mathbf{m}^\psi_{(\alpha \beta)} + \delta \hat{\mathbf{u}}^T b^\lambda_{(\alpha} \mathbf{m}^\varpi_{\beta) \lambda},  \\
\delta \varphi_\alpha &= \delta \hat{\boldsymbol{\psi}}^T \mathbf{N}_\psi \mathbf{t}_\alpha +  \delta \hat{\mathbf{u}}^T \mathbf{N}_{u,\alpha} \mathbf{n} .
\end{align}
From these expressions, we obtain the tensorial form of the internal force vectors:
\begin{align}
\mathbf{f}^u_{int} &= \left( \mathbf{m}^u_{(\alpha \beta)} - b_{\alpha \beta} \mathbf{N}_u \mathbf{n} \right) n^{\alpha \beta} + b^\lambda_{(\alpha} \mathbf{m}^\varpi_{\beta) \lambda} m^{\alpha \beta} + \mathbf{N}_{u,\alpha} \mathbf{n} q^\alpha , \label{eq:internal_force_tensor_u} \\
\mathbf{f}^\psi_{int} &= -\mathbf{m}^\psi_{(\alpha \beta)} m^{\alpha \beta} + \mathbf{N}_\psi \mathbf{t}_\alpha q^\alpha .
\label{eq:internal_force_tensor_psi}
\end{align}

Applying the $\mathcal{M}$-operator to the preceding tensorial expressions yields the explicit B-matrices required for the final matrix-vector form: 
\begin{align}
\mathbf{B}^n &= \mathcal{M} ( \mathbf{m}^u_{(\alpha \beta)} - b_{\alpha \beta} \mathbf{N}_u \mathbf{n} ), \\
\mathbf{B}^m_u &= \mathcal{M} (b^\lambda_{(\alpha} \mathbf{m}^\varpi_{\beta) \lambda}) \, , \quad \mathbf{B}^m_\psi = \mathcal{M} (\mathbf{m}^\psi_{(\alpha \beta)}), \\
\mathbf{B}^q_u &= \mathbf{N}^q \, , \quad \mathbf{B}^q_\psi = \mathbf{N}_\psi \mathbf{T}, \\
\mathbf{N}^q &= \begin{bmatrix} \mathbf{N}_{u,1} \mathbf{n} & \mathbf{N}_{u,2} \mathbf{n} \end{bmatrix} \, , \, \mathbf{T} = \begin{bmatrix} \mathbf{t}_1 & \mathbf{t}_2 \end{bmatrix}.
\end{align}
Ultimately, the internal forces can be expressed compactly as:
\begin{align}
\mathbf{f}_{int}^u &= \mathbf{B}^n \mathcal{M} (n^{\alpha \beta}) + \mathbf{B}^m_u \mathcal{M} (m^{\alpha \beta}) + \mathbf{B}^q_u \mathbf{q} \label{internalf1}\\
\mathbf{f}_{int}^\psi &= -\mathbf{B}^m_\psi \mathcal{M} (m^{\alpha \beta}) + \mathbf{B}^q_\psi \mathbf{q}. \label{internalf2}
\end{align}

Next, we proceed with the linearization to derive the tangent stiffness matrix. Since the B-matrices are independent of displacement and rotation under our assumptions, the linearization of the internal forces in Eqs.~\eqref{eq:internal_force_tensor_u}--\eqref{eq:internal_force_tensor_psi} is straightforward:
\begin{align}
\Delta \mathbf{f}^u_{int} &= \left( \mathbf{m}^u_{(\alpha \beta)} - b_{\alpha \beta} \mathbf{N}_u \mathbf{n} \right) \Delta n^{\alpha \beta} \nonumber + b^\lambda_{(\alpha} \mathbf{m}^\varpi_{\beta) \lambda} \Delta m^{\alpha \beta} + \mathbf{N}_{u,\alpha} \mathbf{n} \Delta q^\alpha, \\
\Delta \mathbf{f}^\psi_{int} &= -\mathbf{m}^\psi_{(\alpha \beta)} \Delta m^{\alpha \beta} + \mathbf{N}_\psi \mathbf{t}_\alpha \Delta q^\alpha.
\end{align}
This requires the linearization of the membrane forces, bending moments, and shear forces, which are obtained from Eqs.\eqref{mforce}--\eqref{qforce} as follows:
\begin{align}
\Delta n^{\alpha \beta }&=2(\sigma \Delta \gamma^\lambda _\lambda a^{\alpha \beta } + \Delta \gamma^{\alpha \beta})-\frac{1}{3}\Bigl[ \Delta \rho^{(\alpha \lambda }b^{\prime \beta )}_\lambda + \sigma a^{\alpha \beta }\Bigl( b_{\mu \nu} \Delta\rho^{\mu \nu} \notag
\\
&+\Bigl(\frac{6}{5}\sigma -1\Bigr)H \Delta\rho^\lambda_\lambda \Bigr)+\frac{3}{5}\sigma \Delta \rho^\lambda _\lambda b^{\alpha \beta }\Bigr], \label{dmforce}
\\
\Delta m^{\alpha \beta }&=\frac{1}{6}(\sigma \Delta \rho^\lambda _\lambda a^{\alpha \beta }+\Delta \rho^{\alpha \beta})-\frac{1}{3}\Bigl[ \Delta \gamma^{(\alpha \lambda }b^{\prime \beta )}_\lambda+\sigma a^{\alpha \beta }\Bigl( \frac{3}{5}b_{\mu \nu} \Delta \gamma^{\mu \nu} \notag
\\
&+\Bigl(\frac{6}{5}\sigma -1\Bigr)H \Delta \gamma^\lambda_\lambda \Bigr)+\sigma \Delta \gamma^\lambda _\lambda b^{\alpha \beta }\Bigr], \label{dmoments}
\\
\Delta q^\alpha &=\frac{5}{6}a^{\alpha \beta }(\Delta u_{,\beta }+b^\lambda _\beta \Delta u_\lambda + \Delta \psi _\beta ). \label{dqforce}
\end{align}

The remaining task is to express the linearization of the strain measures ($\Delta \gamma^{\alpha \beta}$, $\Delta \rho^{\alpha \beta}$ etc.) in terms of the linearized primary unknowns, i.e., $\Delta \vb{u}$ and $\Delta \vb*{\psi}$ and their derivatives. These relations are listed as follows:
\begin{align}
\Delta \gamma^{\alpha \beta} &= a^{\alpha \lambda} a^{\beta \mu} \Delta \gamma_{\lambda \mu} \, , \, \Delta \gamma^{\alpha}_{\beta} = a^{\alpha \lambda} \Delta \gamma_{\lambda \beta} \, , \, \Delta \gamma_{\alpha \beta} = \Delta u_{(\alpha;\beta)} - b_{\alpha \beta} \Delta u, \\
\Delta \rho^{\alpha \beta} &= a^{\alpha \lambda} a^{\beta \mu} \Delta \rho_{\lambda \mu} \, , \, \Delta \rho^{\alpha}_{\beta} = a^{\alpha \lambda} \Delta \rho_{\lambda \beta} \, , \, \Delta \rho_{\alpha \beta} = -\Delta \psi_{(\alpha;\beta)} + b^\lambda_{(\alpha} \Delta \varpi_{\beta)\lambda}, \\
\Delta \varphi_\alpha &= \Delta \psi_\alpha + \mathbf{n} \cdot \Delta \mathbf{u}_{,\alpha},
\end{align}
in which
\begin{align}
&\Delta u_{\alpha;\beta} = \Delta u_{\alpha,\beta} - \Gamma^\lambda_{\alpha \beta} \Delta u_\lambda \, , \quad \Delta \psi_{\alpha;\beta} = \Delta \psi_{\alpha,\beta} - \Gamma^\lambda_{\alpha \beta} \Delta \psi_\lambda , \label{eq:dvarupsi}\\
&\Delta \varpi_{\alpha \beta} = \dfrac{1}{2} \left( \Delta u_{\beta,\alpha} - \Delta u_{\alpha,\beta} \right), \\
&\Delta \psi_\alpha = \mathbf{t}_\alpha \cdot \Delta \boldsymbol{\psi} \, , \quad \Delta u_\alpha = \mathbf{t}_\alpha \cdot \Delta \mathbf{u} \, , \quad \Delta u = \mathbf{n} \cdot \Delta \mathbf{u}
\label{eq:dvar}
\end{align}
and
\begin{equation}
\Delta u_{\alpha,\beta} = \mathbf{t}_{\alpha,\beta} \cdot \Delta \mathbf{u} + \mathbf{t}_\alpha \cdot \Delta \mathbf{u}_{,\beta} \, , \quad \Delta \psi_{\alpha,\beta} = \mathbf{t}_{\alpha,\beta} \cdot \Delta \boldsymbol{\psi} + \mathbf{t}_\alpha \cdot \Delta \boldsymbol{\psi}_{,\beta}
\label{eq:dvar1}
\end{equation}

The terms in Eqs.~\eqref{eq:dvarupsi}--\eqref{eq:dvar1} are ultimately functions of $(\Delta \mathbf{u},\Delta \boldsymbol{\psi})$ and the corresponding spatial derivatives. These quantities are, in turn, interpolated using the shape function matrices (e.g., $\Delta \mathbf{u} = \mathbf{N}_u \Delta \hat{\mathbf{u}}$). By substituting these interpolations back through the preceding equations, the linearization of the internal forces is completed. Finally, the full tangent stiffness matrix is assembled in its matrix-vector form by applying the same mapping operator, $\mathcal{M}$, used for the internal force vector. Its explicit expression is therefore omitted for brevity.

\end{document}